\documentclass[12pt]{article}
\pdfoutput=1
\usepackage[top=1in, bottom=1in, left=0.9in, right=0.9in]{geometry}
\usepackage{amsmath,amsfonts,amssymb,amsthm}
\usepackage{hyperref}
\usepackage{mathrsfs,dsfont}
\usepackage{graphicx,color}
\usepackage{colortbl,dcolumn}
\usepackage{psfrag}
\usepackage{booktabs}
\usepackage{cite}
\usepackage{url}
\usepackage{enumerate}
\allowdisplaybreaks[4]
\numberwithin{equation}{section}
\usepackage[justification=centering]{caption}
\usepackage{authblk}
\usepackage{booktabs}
\usepackage{threeparttable}
\usepackage{cleveref}
\usepackage{graphicx}
\usepackage{subcaption}

\newtheorem{theorem}{Theorem}%
\newtheorem{proposition}[theorem]{Proposition}
\newtheorem{definition}[theorem]{Definition}
\newtheorem{lemma}[theorem]{Lemma}

\newtheorem{remark}{Remark}%

\numberwithin{equation}{section}

\begin{document}

\title{A Positivity-Preserving Numerical Method for the Free Cox-Ingersoll-Ross Equation}

\author[a]{Yuanling Niu}
\author[b]{Zixin Peng}
\author[a]{Jiaxin Wei}
\author[a]{Zhi Yin}
\author[a]{Dan Zeng \thanks{Corresponding author: zengdan@csu.edu.cn}}
\affil[a]{School of Mathematics and Statistics, HNP-LAMA, Central South University, Changsha, 410083, China.}
\affil[b]{School of Mathematics and Statistics, Chongqing University, Chongqing, 400044, China.}
\renewcommand*{\Affilfont}{\small \it}
\date{}

\maketitle

\begin{abstract}
We introduce a positivity-preserving drift-implicit square-root  Euler method for the free Cox--Ingersoll--Ross (CIR) equation driven by free Brownian motion. The method is obtained by discretizing the square-root equation implicitly in the drift and preserves strict positivity of the numerical solution without any step-size restriction. Under the Feller condition, we establish the uniqueness of the global strictly positive solution and  derive an upper bound for its operator norm. 
A key ingredient for the error analysis is a free stochastic differential representation of the inverse square-root process.
Another key ingredient is an \textit{a priori} grid-point error estimate for sufficiently small step sizes that yields uniform upper and strictly positive lower spectral bounds for the numerical solution on any fixed finite time interval. These results support the strong convergence of order $1$ and improve the classical results.
Moreover, both the exact and numerical solutions are exponentially stable with respect to perturbations of the initial value. The numerical stability holds   for arbitrary step
sizes, and the associated discrete decay rate converges to its continuous counterpart as the step size tends to zero. Numerical experiments confirm the theoretical results.
\end{abstract}


\section{Introduction}
This paper  proposes and investigates a structure-preserving numerical scheme for the free CIR equation, with positivity and consequently self-adjointness preserved at the discrete level. The equation takes the form
\begin{align}\label{eq:freeCIRmain}
\mathrm{d}X_t=(a-bX_t)\mathrm{d}t+\frac{\sigma}{2}\sqrt{X_t}\mathrm{d}W_t+\frac{\sigma}{2}\mathrm{d}W_t\sqrt{X_t}, \quad X_0\in \mathcal{A}_+,
\end{align}
where $\mathcal{A}_+$ denotes the strictly positive and self-adjoint elements of a finite von Neumann algebra $\mathcal{A}$ equipped with a faithful normal tracial state $\varphi$,  $a, b, \sigma >0$ are scalar constants and $\{W_t\}_{t\geq0}$ is a free Brownian motion.
When the Feller condition \(2a\ge \sigma^2\) holds,  Graf, Port and Schl\"uchtermann \cite{freecir}  established the global existence of a strictly positive solution. 

The free CIR equation is an important type of free stochastic differential equations (SDEs). They describe the evolution of operators in free probability, where classical independence is replaced by free independence \cite{freeprobability,Voicu2}. The theory of free SDEs was developed through the free stochastic calculus of Biane-Speicher \cite{fff,Stochasticcalculus} and subsequently in the work of Kargin \cite{onfree}. One important motivation for studying free SDEs comes from random matrix theory: free SDEs can be regarded as limits of matrix-valued SDEs. For instance, suitably normalized Hermitian matrix-valued Brownian motions converge to free Brownian motion \cite{Voicu}, while related random-matrix limits include free Wishart processes \cite{freewishart}, free diffusions with convex potential \cite{GuionnetShlyakhtenko2009}, and general free SDEs with suitable assumptions \cite{matrixtofree}. In this sense, the free CIR equation can be viewed as the limit of matrix-valued CIR processes \cite{matrixtofree}, which provides potential applications to the spectral analysis of high-dimensional stochastic covariance models in multivariate finance.

The classical CIR equation is widely used in models of interest rates, stochastic volatility, and credit risk \cite{cox1985theory,heston1993closed,duffee1999estimating}. There is a sizable literature dedicated to CIR-type numerical approximations within the classical framework (e.g., \cite{cir1,Alfonsi2010,Alfonsi,cirnum2,d-i,cir4,cir5,HefterJentzen2019,hong,MoroSchurz2007,cirnum}), with various schemes constructed specifically for its square-root diffusion coefficient. In particular, 
Dereich, Neuenkirch and Szpruch \cite{d-i} established strong convergence of order $1/2$ at the grid points under the strict Feller condition \(2a>\sigma^2\). Alfonsi \cite{Alfonsi}, Neuenkirch and Szpruch \cite{cirnum} subsequently established strong convergence of order 1 under stronger parameter restrictions than the strict Feller condition. Recently, following this line of research, Liu, Neuenkirch and Wang \cite{LiuNeuenkirchWang2025} improved the convergence order to $3/2$ under additional assumptions. 
This naturally raises the question of whether one can develop a numerical method in the free setting that preserves the strict positivity of the solution.

However, due to the noncommutativity of the operators, the classical methods cannot be directly adapted to the free setting. To the best of our knowledge, only two works \cite{numfsde1,numfsde2} have developed numerical schemes for free SDEs with locally operator Lipschitz coefficients. In \cite{numfsde1},  the authors established strong convergence of order $1/2$ in $L^1(\varphi)$ for the Euler–Maruyama scheme, while the authors of \cite{numfsde2} established strong convergence of order $1/2$ in $L^2(\varphi)$ for the $\theta$-method. However, the free CIR equation lies outside those frameworks, and their numerical schemes, when applied to the free CIR equation, do not in general preserve the positivity.
In this paper, we construct a drift-implicit Euler method. Under the Feller condition, we prove strong convergence of the method and establish long-time exponential stability in $L^2(\varphi)$ with respect to perturbations of the initial value for the exact and numerical solutions. More importantly, the proposed method preserves strict positivity of the solution; see \Cref{prop:wellposed}. Beyond the free CIR setting, the numerical construction and analytical techniques derived here may be capable to further studies of structure-preserving numerical methods for free SDEs with non-Lipschitz coefficients.

Our strategy is to introduce the square-root process $V_t = \sqrt{X_t}$, which makes the noise additive but introduces \(V_t^{-1}\) into the drift. 
A simple observation is that $X_t$ admits a uniform positive lower spectral bound on every fixed finite time interval (this bound depends on the solution and on the time interval; see \Cref{lem:spectral_gap}).
Therefore, on the subset of strictly positive operators, the square-root map is locally operator Lipschitz. By using free stochastic calculus, we show that $V_t^{-1}$ is the solution of the free SDE given in \eqref{eq:inverse}; see Theorem \ref{prop:inverse}. 
Based on this result, our analysis is divided to the following steps: {\it Step 1.} We decompose the remainder term $r_k$ (given in \eqref{eq:rk}) into finite-variation and martingale components. It allows the use of the conditional cancellation of martingale increments in $L^2(\varphi)$, which is essential to improve the error estimate. {\it Step 2.}  We compare the exact and numerical square-root processes at temporal grid points. We first prove that, for sufficiently small step sizes, an operator norm error estimate of order $1/2$ holds, and this estimate transfers the exact solution's spectral bounds to the numerical approximation on any fixed finite interval, where these bounds are independent of the step size. {\it Step 3.} Combining the decomposition of  \(r_k\) with the numerical upper bound gives order 1 strong convergence in $L^2(\varphi)$; see \Cref{thm:convergence}. We remark that our results hold with the Feller condition, and thus improve the classical results in \cite{d-i,Alfonsi,cirnum, LiuNeuenkirchWang2025}. 
Moreover, we establish exponential stability in $L^2(\varphi)$ for exact and numerical solutions with perturbed initial values, which holds for arbitrary step sizes, and whose decay rate converges to the continuous rate as the step size tends to zero; see \Cref{thm:stability1,thm:numstability}. Numerical experiments illustrate the above convergence rates, stability, and preservation of positivity.

The rest of the paper is organized as follows: 
 \Cref{sec:pre} presents the preliminaries. \Cref{sec:main} provides the analytical estimates for the free CIR process, and \Cref{sec:method} develops the drift-implicit square-root Euler scheme and proves that it preserves positivity.  \Cref{sec:results} proves convergence and stability, and  \Cref{sec:experiments} reports the numerical experiments.

\section{Preliminaries}
\label{sec:pre}

This section briefly reviews the required concepts from free probability and free stochastic calculus. 

\subsection{Free probability}
Let $(\mathcal{A}, \{\mathcal{A}_t\}_{t\geq0}, \varphi)$ be a noncommutative probability space, where $\mathcal{A}$ is a finite von Neumann algebra, $\{\mathcal{A}_t\}_{t\geq0}$ is a filtration of subalgebras of $\mathcal{A}$, i.e.,  an increasing family with $\mathcal{A}_s\subset\mathcal{A}_t$ for $s\leq t$ and $\varphi\colon \mathcal{A} \to \mathbb{C}$ is a faithful, normal, and tracial state.

An operator valued process
$\{X_t\}_{t\geq0}$ is adapted if $X_t\in\mathcal{A}_t$ for every $t$. Following \cite{Non-C}, for a von Neumann subalgebra $\mathcal{B}\subseteq\mathcal{A}$, the $\varphi$-preserving conditional expectation $\mathbb{E}_{\mathcal{B}}:\mathcal{A}\to\mathcal{B}$
satisfies $\mathbb{E}_{\mathcal{B}}(b)=b, \mathbb{E}_{\mathcal{B}}(b_1xb_2)=b_1\mathbb{E}_{\mathcal{B}}(x)b_2,$ for all $x\in\mathcal{A}$ and $b,b_1,b_2\in\mathcal{B}$, and
$\varphi\bigl(\mathbb{E}_{\mathcal{B}}(x)\bigr)=\varphi(x).$
For the filtration $\{\mathcal{A}_t\}_{t\geq0}$, we write
$\mathbb{E}_t:=\mathbb{E}_{\mathcal{A}_t}$. 

The subalgebras $\mathcal{A}_1, \dots, \mathcal{A}_m$ of $\mathcal{A}$ are said to be freely independent if 
$\varphi(a_{i(1)} \cdots a_{i(k)}) = 0$ whenever $k\geq 1$, $i(1),\ldots,i(k)\in\{1,\ldots,m\}$, $a_{i(j)}\in\mathcal{A}_{i(j)}$ with $\varphi(a_{i(j)})=0$ for all $j=1,\ldots,k$, and $i(j)\neq i(j+1)$ for all $j=1,\ldots,k-1$ \cite[Definition 11]{Non-C}.

 We write $\mathcal{A}_{sa}$ for the set of self-adjoint elements of $\mathcal{A}$. For $x\in\mathcal{A}_{sa}$, $\operatorname{spec}(x)$ denote the set of $\lambda\in\mathbb{C}$ for which $\lambda-x$ fails to be invertible in $\mathcal{A}$; its complement $\rho(x)=\mathbb{C}\setminus\operatorname{spec}(x)$ is the resolvent set; $\mathcal{A}_+$ denotes the set of $\mathcal{A}_{sa}$  with strictly positive spectrum. For $\lambda\in\rho(x)$, the resolvent is defined as $R(x,\lambda)=(\lambda-x)^{-1}$. The spectral distribution of $x$ is the unique probability measure $\mu$ on $\mathbb{R}$ satisfying $\varphi(x^k)=\int_{\mathbb{R}}t^k\,\mathrm{d}\mu(t)$ for all $k\in\mathbb{N}$; refer to \cite{Non-C}.

For $1 \leq p < \infty$, the completion of $\mathcal{A}$ with respect to the norm $\|x\|_p := [\varphi(|x|^p)]^{^{\frac{1}{p}}}$ yields the noncommutative $L^p(\varphi)$. The $L^\infty(\varphi)$ is defined as $\mathcal{A}$ itself, equipped with the operator norm $\|\cdot\|$. Since $\varphi$ is unital, the norms are monotone in $p$: $\|x\|_p \leq \|x\|_q,\ 1 \leq p \leq q \leq \infty,$
so in particular $\|x\|_2 \leq \|x\|$.

We shall use the following standard inequalities in noncommutative $L^p(\varphi)$; see \cite{Non-C,freecheb}. The noncommutative H\"older inequality asserts that for any $0 < r, p,q \leq \infty$ with $\frac{1}{r}=\frac{1}{p} + \frac{1}{q},$ if $x\in L^p(\varphi)$, $y\in L^q(\varphi)$, then $xy\in L^r(\varphi)$ and
\[
\|xy\|_r\leq\|x\|_p\|y\|_q,\quad \text{i.e.,}\quad \varphi(|xy|^r)^{\frac1r}\leq(\varphi(|x|^p))^{\frac{1}{p}}(\varphi(|y|^q))^{\frac{1}{q}}.
\]
Additionally, we have the noncommutative Chebyshev-type inequality: for any $x \in \mathcal{A}_{sa}$  with spectrum in an interval $J$ and for any real-valued monotone functions $f, g$ on $J$ with the same monotonicity, $\varphi(f(x)g(x))\geq \varphi(f(x))\varphi(g(x)).$

\subsection{Free calculus and free SDEs} 

\begin{definition}
Free Brownian motion is a family $\{W_t\}_{t \geq 0}$ of self-adjoint elements of \(\mathcal{A}\) satisfying the following conditions:
\begin{enumerate}
    \item $W_0 = 0$;
    \item for $t > s \geq 0$, the increment $W_t - W_s$ is freely independent  of the subalgebra $\mathcal{A}_s$;
    \item $W_t - W_s$ has a semicircular distribution with mean $0$ and variance $t-s$.
\end{enumerate}
\end{definition}

Fix $T>0$, assume that $\beta_t,\gamma_t \in \mathcal{A}_t$ for all $t\in[0,T]$ and the function $t \mapsto \|\beta_t\|\cdot\|\gamma_t\|$ belongs to $L^2([0,T])$, the It\^o-style free stochastic integral \( \int_0^T \beta_t \, \mathrm{d}W_t \, \gamma_t\) is well-defined; see \cite{Stochasticcalculus, onfree} for further details here and below.

As the free probability counterpart of the classical It\^o formulas, for adapted processes $a_t,b_t,c_t,e_t,$ we have the following formulas:
\begin{align*}
(a_t\,\mathrm{d}t)(b_t\,\mathrm{d}t) &= 0, & (a_t\,\mathrm{d}t)(b_t\,\mathrm{d}W_t\,c_t) &= 0, \\
(a_t\,\mathrm{d}W_t\,b_t)(c_t\,\mathrm{d}t) &= 0, & (a_t\,\mathrm{d}W_t\,b_t)(c_t\,\mathrm{d}W_t\,e_t) &= \varphi(b_tc_t)\,a_te_t\,\mathrm{d}t.
\end{align*}
The free stochastic integral satisfies a free analogue of the Burkholder--Gundy (BG) inequality  and free It\^o isometry:
\[
\bigg\|\int_0^Tb_t\,\mathrm{d}W_t\,c_t\bigg\|^2 \leq 8\int_0^T\|b_t\|^2\cdot\|c_t\|^2\,\mathrm{d}t;\quad\bigg\|\int_0^Tb_t\mathrm{d}W_tc_t\bigg\|^2_2=\int_0^T\|b_t\|_2^2\cdot\|c_t\|_2^2\,\mathrm{d}t.
\]

\begin{definition}
A general class of free SDEs is given by
\begin{align}\label{eq:freeSDE}
\mathrm{d}X_t = \alpha(X_t)\mathrm{d}t + \sum_{i=1}^m \beta_i(X_t)\mathrm{d}W_t\gamma_i(X_t),
\end{align}
where the coefficients $\alpha, \beta_i, \gamma_i:\mathcal{A}_{sa} \to \mathcal{A}_{sa}$ are continuous functions in the operator norm.
A solution to \eqref{eq:freeSDE} with initial condition $X_0$ is an adapted 
process $\{X_t\}_{0\le t\le T}$ satisfying the following properties:
\begin{enumerate}
\item $X_t \in \mathcal{A}_{sa}$ for all $0 \leq t \leq T$;
\item $\|\alpha(X_t)\| \in L^1([0,T])$ and $\|\beta_i(X_t)\|\cdot\|\gamma_i(X_t)\| \in L^2([0,T]), i = 1,\dots,m$;
\item $X_t$ satisfies the integral equation for all $0 \leq t \leq T$:
\begin{align}\label{eq:freeInt}
X_t = X_0 + \int_0^t \alpha(X_s)\mathrm{d}s + \sum_{i=1}^m \int_0^t \beta_i(X_s)\mathrm{d}W_s\gamma_i(X_s),
\end{align}
\noindent\hspace*{-3.5em}where the stochastic integrals are understood in terms of the free It\^o integral.
\end{enumerate}
\end{definition}

In this paper, we will also use the functional free It\^o formula for the solution $X_t.$ Recall the noncommutative derivations as follows \cite{Stochasticcalculus}: 
Let \(\mathcal A^{\mathrm{op}}\) denote the opposite algebra of \(\mathcal A\).
The algebraic tensor product
\(\mathcal A\otimes\mathcal A^{\mathrm{op}}\) acts on \(\mathcal A\) and \(\mathcal A\otimes\mathcal A^{\mathrm{op}}\), respectively, via
\begin{align}\label{eq:sharp}
    \begin{aligned}
        (a\otimes b)\sharp u:=aub,\qquad (a\otimes b)\sharp(u\otimes v):=au\otimes vb,\qquad a,b,u,v\in\mathcal A.
    \end{aligned}
\end{align}
Let \(\partial:\mathbb C[X]\longrightarrow\mathbb C[X]\otimes\mathbb C[X]\) be the canonical derivation, defined on monomials by
\[
\partial X^n=\sum_{k=0}^{n-1}
X^k\otimes X^{\,n-k-1},
\qquad n\ge1.
\]
As in \cite{Voiculescu2000Coalgebra}, the second derivative \(\partial^2:\mathbb C[X]\longrightarrow\mathbb C[X]\otimes\mathbb C[X]\otimes\mathbb C[X]\)
is defined by:
\begin{equation*}
\begin{split}
\partial^2X^n &:=2 (\partial \otimes \mathrm{id}) \circ \partial X^n = 2\sum_{\substack{k,l\ge0\\k+l\le n-2}}
X^k\otimes X^l\otimes X^{\,n-k-l-2},\qquad n\ge2.
\end{split}
\end{equation*}
Let \(m:\mathcal A\otimes\mathcal A\longrightarrow\mathcal A,\)
be the multiplication map, and define the contraction \(\eta:=m\circ\bigl(I_{\mathcal A}\otimes\varphi\otimes I_{\mathcal A}\bigr):\mathcal A\otimes\mathcal A\otimes\mathcal A\longrightarrow\mathcal A.\)
Following \cite[Definition~4.3.1]{Stochasticcalculus}, for
\(X\in\mathcal A\), \(U\in\mathcal A\otimes\mathcal A\), and \(n\geq 0\),
define
\begin{align}\label{eq:deltau1}
    \begin{aligned}
        \Delta_U(X^n):=\eta\left((1\otimes U)\sharp\partial^2(X^n)\sharp(U\otimes1)\right).
    \end{aligned}
\end{align}
Hence, for the process $X_t$ given in \eqref{eq:freeInt}, denote $U_s:=\sum_{i=1}^m  \beta_i(X_s)\otimes \gamma_i(X_s).$  The functional free It\^o formula
\cite[Proposition~4.3.2]{Stochasticcalculus} states that, for every polynomial \(P\in\mathbb C[X]\),
\begin{align}\label{eq:polyfreeito}
    \begin{aligned}
    P(X_t)=P(X_0)+\int_0^t \partial P(X_s)\sharp \mathrm{d}X_s+\frac12\int_0^t \Delta_{U_s}P(X_s)\mathrm{d}s.
    \end{aligned}
\end{align}
\section{Free CIR equation}
\label{sec:main}
In the sequel, we assume that the coefficients of the free CIR equation \eqref{eq:freeCIRmain} satisfy:
\[
a,b,\sigma>0,\quad 2a\geq \sigma^2,\quad X_0\in\mathcal A_+.
\]
Under these assumptions, the existence of a global strictly positive solution to \eqref{eq:freeCIRmain} follows from \cite{freecir}. Hence, there exists a solution $\{X_t\}_{t\geq 0}$ satisfying \(X_t\in\mathcal{A}_+,\ t\geq 0.\)
Moreover, the solution is continuous in \(t\) with respect to the operator norm.

Now, we establish a uniform operator norm bound for the solution, which has not been studied in prior work and is useful for the subsequent analysis.
\begin{proposition}
\label{prop:uniform_bound}
Let $\{X_t\}_{t \geq 0}$ be a global positive solution of the free CIR equation
\eqref{eq:freeCIRmain}. Then there exists a constant $M>0$, independent of $t$, such that
\[
\sup_{t \geq 0}\|X_t\|\leq M,\ \text{where}\ M=\left(
\frac{\sigma}{\sqrt{b}}
+
\sqrt{
\|X_0\|+\frac{\sigma^2+a}{b}}
\right)^2.
\]

\end{proposition}
\begin{proof}
Applying the free It\^{o} formula to $e^{bt}X_t$, we obtain
\[
\begin{aligned}
\mathrm{d}(e^{bt}X_t)
&=
be^{bt}X_t\,\mathrm{d}t
+
e^{bt}\,\mathrm{d}X_t\\
&=
ae^{bt}\,\mathrm{d}t
+
\frac{\sigma}{2}e^{bt}\sqrt{X_t}\,\mathrm{d}W_t
+
\frac{\sigma}{2}e^{bt}\mathrm{d}W_t\sqrt{X_t}.
\end{aligned}
\]
Integrating both sides over $[0,t]$ and then multiplying by $e^{-bt}$, we obtain
\[
\begin{aligned}
X_t
={}&
e^{-bt}X_0
+
\frac{a}{b}(1-e^{-bt})+
\frac{\sigma}{2}e^{-bt}
\int_0^t e^{bs}\sqrt{X_s}\,\mathrm{d}W_s+\frac{\sigma}{2}e^{-bt}
\int_0^t e^{bs}\mathrm{d}W_s\sqrt{X_s}.
\end{aligned}
\]
Therefore,
\[
\begin{aligned}
\|X_t\|
\leq{}&
\|X_0\|
+
\frac{a}{b}+
\frac{\sigma}{2}e^{-bt}
\bigg(\left\|
\int_0^t e^{bs}\sqrt{X_s}\,\mathrm{d}W_s
\right\|+
\left\|
\int_0^t e^{bs}\mathrm{d}W_s\sqrt{X_s}
\right\|\bigg).
\end{aligned}
\]
By the free BG inequality in the operator norm, we have
\[
\left\|
\int_0^t e^{bs}\sqrt{X_s}\,\mathrm{d}W_s
\right\|+\left\|
\int_0^t e^{bs}\mathrm{d}W_s\sqrt{X_s}
\right\|
\leq
4\sqrt{2}
\left(
\int_0^t e^{2bs}\|\sqrt{X_s}\|^2\,\mathrm{d}s
\right)^{\frac12}.
\]
Since $X_s$ is positive,
$
\|\sqrt{X_s}\|^2=\|X_s\|.
$
Consequently,
\[
\|X_t\|
\leq
\|X_0\|
+
\frac{a}{b}
+
2\sqrt{2}\sigma e^{-bt}
\left(
\int_0^t e^{2bs}\|X_s\|\,\mathrm{d}s
\right)^{\frac12}.
\]
For $t\geq0$, set
$
M_t:=\sup_{0\leq s\leq t}\|X_s\|.
$
Since $s\mapsto X_s$ is continuous in the operator norm, $M_t<\infty$
for every finite $t$. Moreover,
\[
\int_0^t e^{2bs}\|X_s\|\,\mathrm{d}s
\leq
M_t\int_0^t e^{2bs}\,\mathrm{d}s
=
M_t\frac{e^{2bt}-1}{2b}.
\]
It follows that
\[
\|X_t\|
\leq
\bigg(\|X_0\|+\frac{a}{b}\bigg)+\frac{2\sigma}{\sqrt{b}}\sqrt{M_t}.
\]
Note that $M_r\leq M_t,\ 0\leq r\leq t$. 
Taking the supremum over $r\in[0,t]$ gives
\[
M_t
\leq
\bigg(\|X_0\|+\frac{a}{b}\bigg)+\frac{2\sigma}{\sqrt{b}}\sqrt{M_t}.
\]
Setting $y_t:=\sqrt{M_t}$, we obtain
$y_t^2-\frac{2\sigma}{\sqrt{b}}y_t-\bigl(\|X_0\|+\frac{a}{b}\bigl)\leq0.$
Since $y_t\geq0$, it follows that
\[
y_t
\leq
\frac{\sigma}{\sqrt{b}}
+
\sqrt{\|X_0\|+\frac{\sigma^2+a}{b}}.
\]
Therefore,
since the right-hand side is independent of $t$,
we conclude that
\[
\sup_{t\geq0}\|X_t\|
\leq
\left(
\frac{\sigma}{\sqrt{b}}
+
\sqrt{\|X_0\|+\frac{\sigma^2+a}{b}}\right)^2
<\infty.
\]
This completes the proof.
\end{proof}

\subsection{Uniqueness of free CIR equation} 
In 2011, Kargin \cite{onfree} established local well-posedness for the free SDEs \eqref{eq:freeSDE} under the assumption that the coefficients are locally operator-Lipschitz. However, the free CIR equation \eqref{eq:freeCIRmain} does not satisfy this condition due to its square-root diffusion coefficient, and therefore falls outside Kargin's scope.
In 2022, Garf, Port and Schl\"uchterma \cite{freecir} obtained a global existence result for the free CIR equation \eqref{eq:freeCIRmain}. In this subsection, we aim to prove uniqueness of the global strictly positive solution to \eqref{eq:freeCIRmain}.

The proof of uniqueness relies on the following two preliminary lemmas. The first one is that, on every finite closed time interval, the strict positivity and operator norm continuity of the solution yield a uniform spectral separation from zero, with the separation constant depending on the solution and the time interval .
\begin{lemma}
\label{lem:spectral_gap}
Let $\{X_t\}_{t\geq0}$ be a global solution of the free CIR equation \eqref{eq:freeCIRmain}. Then, for every $T>0$, there exist two constants $0<m^X_T\leq M^X$ such that
\[
m^X_T\mathbf{1}\leq X_t\leq M^X\mathbf{1},
\quad 0\leq t\leq T.
\]
Equivalently,
\[
0<\inf_{0\leq t\leq T}\inf\operatorname{spec}(X_t)\leq\sup_{0\leq t\leq T}\sup\operatorname{spec}(X_t)<\infty.
\]
\end{lemma}
\begin{proof} According to \(
-\|A\|\mathbf{1}\leq A\leq \|A\|\mathbf{1}, \ \forall{A\in\mathcal{A}_{sa}},
\)
and $t\mapsto X_t$ is continuous in the operator norm,
it is  easy to obtain the map
$
t\mapsto\inf\operatorname{spec}(X_t)$ is continuous on $[0,T].$ 
Combining this estimate with \Cref{prop:uniform_bound}, the remaining proof follows directly, and we omit the details. 
\end{proof}

The lemma below is an immediate consequence of the van Hemmen--Ando inequality established in \cite{vanhemmen}; we omit the proof. 
\begin{lemma}
\label{lem:sqrt_lipschitz}
Let $A,B\in\mathcal{A}_+$ satisfy
$
A,B\geq m\mathbf{1},\  \forall m>0.
$
Then
\[
\|\sqrt{A}-\sqrt{B}\|
\leq
\frac{1}{2\sqrt{m}}\|A-B\|.
\]
\end{lemma}

\begin{proposition}
    \label{thm:uniqueness}
Let $X=\{X_t\}_{t\geq0}$ and $Y=\{Y_t\}_{t\geq0}$ be two global strictly positive solutions of the free CIR equation \eqref{eq:freeCIRmain} with the same initial value $X_0=Y_0$, then
\[
X_t=Y_t,\quad t\geq0.
\]
\end{proposition}
\begin{proof}
For any $T>0$, according to \Cref{lem:spectral_gap}, there exist constants
$m_T^X>0$ and $m_T^Y>0$ such that
\[
X_t\geq m_T^X\mathbf{1},
\quad
Y_t\geq m_T^Y\mathbf{1},
\quad 0\leq t\leq T.
\]
Set
$
m_T:=\min\{m_T^X,m_T^Y\}>0,
$
and \Cref{lem:sqrt_lipschitz} yields
\[
\|\sqrt{X_t}-\sqrt{Y_t}\|
\leq
\frac{1}{2\sqrt{m_T}}\|X_t-Y_t\|,
\quad 0\leq t\leq T.
\]
For simplicity, denote $
L_T:=\frac{1}{2\sqrt{m_T}}$
and $
Z_t:=X_t-Y_t.$
Since $X_0=Y_0$, subtracting the two equations gives
\[
\begin{aligned}
Z_t
={}&
-b\int_0^t Z_s\,\mathrm{d}s+
\frac{\sigma}{2}
\int_0^t
\bigl(\sqrt{X_s}-\sqrt{Y_s}\bigr)\,\mathrm{d}W_s+
\frac{\sigma}{2}
\int_0^t
\mathrm{d}W_s
\bigl(\sqrt{X_s}-\sqrt{Y_s}\bigr).
\end{aligned}
\]
Taking the operator norm, using the triangle inequality and the free BG inequality, we obtain
\[
\begin{aligned}
\|Z_t\|
\leq{}&
b\int_0^t\|Z_s\|\,\mathrm{d}s+
\frac{\sigma}{2}
\left\|
\int_0^t
(\sqrt{X_s}-\sqrt{Y_s})\,\mathrm{d}W_s+\int_0^t
\mathrm{d}W_s(\sqrt{X_s}-\sqrt{Y_s})
\right\|\\
&\leq
b\int_0^t\|Z_s\|\,\mathrm{d}s
+
2\sqrt{2}\sigma
\left(
\int_0^t
\|\sqrt{X_s}-\sqrt{Y_s}\|^2\,\mathrm{d}s
\right)^{\frac12}\\
&\leq b\int_0^t\|Z_s\|\,\mathrm{d}s
+
2\sqrt{2}\sigma L_T
\left(
\int_0^t\|Z_s\|^2\,\mathrm{d}s
\right)^{\frac12}.
\end{aligned}
\]
For $0\leq t\leq T$, using $(u+v)^2\leq2u^2+2v^2$ and the Cauchy--Schwarz inequality, we obtain
\[
\|Z_t\|^2
\leq
\left(
2b^2T+16\sigma^2L_T^2
\right)
\int_0^t\|Z_s\|^2\,\mathrm{d}s.
\]
By Gronwall's inequality, we have
$\|Z_t\|^2=0,\ 0\leq t\leq T.$
Since $T>0$ is arbitrary, we conclude that
$
X_t=Y_t,\
 t\geq0.
$
\end{proof}

\subsection{\texorpdfstring{The square-root process $V$ and its inverse $V^{-1}$}{The square-root process V and its inverse V inverse}} 
To construct a positivity-preserving numerical approximation for \eqref{eq:freeCIRmain}, we introduce an auxiliary square-root process. Define
\[
\xi:=a-\frac{\sigma^2}{4}\quad\ \text{and}\quad f(A):=\frac{\xi}{2}A^{-1}-\frac{b}{2}A,
\quad A\in\mathcal A_+.
\]
Consider the free SDE
\begin{equation}\label{eq:freeCIRlam}
\mathrm{d}V_t=f(V_t)\mathrm{d}t+\frac{\sigma}{2}\mathrm{d}W_t,
\quad
V_0=\sqrt{X_0}.
\end{equation}
By \cite[Theorem~3.3]{freecir},  \eqref{eq:freeCIRlam} admits a strictly positive, operator norm continuous global solution \(\{V_t\}_{t\geq0}\).
Since $\xi+\frac{\sigma^2}{4}=a$, the free It\^o product formula shows that \(\{V_t^2\}_{t\geq0}\) satisfies \eqref{eq:freeCIRmain} with the initial value \(X_0\). By 
\Cref{thm:uniqueness}, we therefore have
\begin{align}\label{eq:square root}
X_t=V_t^2 \quad\text{and}\quad V_t=\sqrt{X_t},
\quad t\geq0.
\end{align}

Hence, the uniform operator norm bound for $X_t$ immediately yields a corresponding bound for $V_t$.
\begin{lemma}\label{lem:V_bound}
     Let \(\{V_t\}_{ t\geq 0}\) be the strictly positive solution of \eqref{eq:freeCIRlam} satisfying $V_t=\sqrt{X_t}$.
    Then, there exists a constant \(M_V>0\), independent of \(t\), such that
    \[
    0<V_t\leq M_V\mathbf 1,\quad t\geq0.
    \]
    Moreover, for any $T>0,$ there exist constants \(\ell_T>0\) and \(K_T<\infty\) such that
    \[
    \ell_T\mathbf 1\leq V_t\leq M_V\mathbf 1,
    \quad
    \|V_t^{-1}\|\leq \ell_T^{-1},
    \quad 0\leq t\leq T,
    \]
    and
    \[
    \|V_t-V_s\|
    \leq K_T |t-s|^{\frac12},
    \quad s,t\in[0,T].
    \]
    In particular, one may take
    $
    M_V=\sqrt{M},\
    \ell_T=\sqrt{m_T},
    $
    where \(M\) is the constant in  \Cref{prop:uniform_bound} and \(m_T=m^X_T\) is the constant in \Cref{lem:spectral_gap} (by the uniqueness, we may drop the superscript and write simply $m_T$).
\end{lemma}

\begin{proof}
The first two assertions follow immediately from $V_t=\sqrt{X_t}$, \Cref{prop:uniform_bound} and \Cref{lem:spectral_gap}. It remains to prove the temporal regularity.
For $0\leq s<t\leq T$, the integral form of~\eqref{eq:freeCIRlam}
gives
\[
    V_t-V_s
    =
    \int_s^t
f(V_r)\mathrm{d}r
    +
    \frac{\sigma}{2}(W_t-W_s).
\]
Hence, using the above bounds and
$\|W_t-W_s\|=2|t-s|^{\frac12}$, we have
\[
\begin{aligned}
    \|V_t-V_s\|
    &\leq
    \left(
        \frac{\xi}{2\ell_T}
        +
        \frac{bM_V}{2}
    \right)|t-s|
    +
    \sigma |t-s|^{\frac12} \leq
    K_T |t-s|^{\frac12}.
\end{aligned}
\]
\end{proof}

The above estimate controls the size of the square-root process.
On the other hand, the inverse term $V_t^{-1}$ appears naturally in the drift coefficient of \eqref{eq:freeCIRlam}.
This, as well as the requirements for what follows, motivates us to further characterize the dynamics of the inverse square-root process.

\begin{theorem}\label{prop:inverse}
Let \(\{V_t\}_{t\geq0}\) be the strictly positive solution of \eqref{eq:freeCIRlam} with initial condition \(V_0\in\mathcal{A}_+\). Then its inverse \(\{V_t^{-1}\}_{t\geq0}\) is a strictly positive solution of the following free SDE
\begin{equation}\label{eq:inverse}
\begin{aligned}
\mathrm{d}Y_t = \Bigl(-\frac{\xi}{2}Y_t^3 &+ \frac{\sigma^2}{4}\varphi(Y_t)(Y_t)^2 + \frac{b}{2}Y_t\Bigr)\mathrm{d}t- \frac{\sigma}{2} Y_t\mathrm{d}W_t Y_t,
\end{aligned}
\end{equation}
with initial condition \(Y_0=V_0^{-1}\in\mathcal{A}_+\).
\end{theorem}

Contrary to the classical case, the free It\^o formula is explicitly stated only for polynomials. Since \(x\mapsto x^{-1}\) is not a polynomial, we cannot directly write down \(\mathrm{d}(V_t^{-1})\) by applying the free It\^o formula. Therefore, we construct a sequence of polynomial approximations of the inverse function via Neumann series expansions on the spectrum of \(V_t\). This approach simultaneously establishes that \(V_t^{-1}\) is a genuine free It\^o process and yields its governing equation.
\begin{proof}[Proof of Theorem \ref{prop:inverse}]
Fix $T>0$. By \Cref{lem:V_bound}, there exist constants \(0 < \ell_T\leq M_V < \infty\) such that
\begin{equation}
\operatorname{spec}(V_t) \subseteq [\ell_T, M_V], \quad 0 \leq t \leq T.
\nonumber
\end{equation}
Define
\[
R_t := \mathbf{1} - \frac{V_t}{ M_V}, \quad \|R_t\| \leq 1 - \frac{\ell_T}{ M_V} =: \rho_T < 1.
\]
The Neumann series for the inverse converges absolutely in operator norm, uniformly in \(t \in [0,T]\):
\begin{equation}\label{eq:vvv}
V_t^{-1} = \frac{1}{M_V}(\mathbf{1} - R_t)^{-1} = \frac{1}{M_V}\sum_{k=0}^{\infty} R_t^{\,k}.
\end{equation}
Define the sequence of partial sums
\[
p_N(x) := \frac{1}{M_V} \sum_{k=0}^{N} \Bigl(1 - \frac{x}{M_V}\Bigr)^k, \quad N \geq 0.
\]
We have
\(p_N(V_t) = \frac{1}{M_V}\sum_{k=0}^{N}R_t^{\,k}\). The geometric tail estimate gives the uniform operator norm convergence: for $N \to \infty,$
\begin{equation}\label{eq:polyconv}
\sup_{0 \leq t \leq T} \bigl\|p_N(V_t) - V_t^{-1}\bigr\|
= \sup_{0 \leq t \leq T} \Bigl\|\frac{1}{M_V}\sum_{k=N+1}^{\infty} R_t^{\,k}\Bigr\|
\leq \frac{1}{M_V}\frac{\rho_T^{N+1}}{1-\rho_T} \longrightarrow 0.
\end{equation}
Since \(V_t\) is adapted to the filtration $\mathcal{A}_t$ and each \(p_N\) is a polynomial, \(p_N(V_t)\) is adapted to the filtration $\mathcal{A}_t$.  Thus, according to  \(\mathcal{A}_t\) is norm-closed and \eqref{eq:polyconv}, we have \(V_t^{-1} \in \mathcal{A}_t, t \in[0,T]\).

By the polynomial free It\^o formula \eqref{eq:polyfreeito} with \eqref{eq:freeCIRlam}, we obtain
\begin{equation}\label{eq:ito-pN}
p_N(V_t)=p_N(V_0)+\int_0^t\partial p_N(V_s)\sharp \bigl(f(V_s)\mathrm ds +\frac{\sigma}{2} \mathrm dW_s\bigl)+\frac12\int_0^t\Delta_{\frac{\sigma}{2}(\mathbf{1} \otimes \mathbf{1})}p_N(V_s)\mathrm ds.
\end{equation} 
We claim that every term in \eqref{eq:ito-pN} converges in operator norm to the corresponding expression of \eqref{eq:inverse}.

For any \(A\) with \(\operatorname{spec}(A)\subset[\ell_T,M_V]\), let \(R_A = \mathbf{1} - A/{M_V}\) such that \(\|R_A\| \leq \rho_T\). 
Then, by the chain rule for the noncommutative polynomial and
\(\partial R_A=-\frac{1}{M_V}(\mathbf{1}\otimes \mathbf{1}),\) we have
\begin{equation}\label{eq:first-der-R}
\partial (R_A^k)=\sum_{j=0}^{k-1} (R_A^{\,j}\otimes\mathbf{1}) \partial R_A (\mathbf{1}\otimes R_A^{\,k-1-j})=-\frac{1}{M_V}  \sum_{j=0}^{k-1} R_A^{\,j} \otimes R_A^{\,k-1-j}.
\end{equation}
Differentiating \(p_N\) termwise, we obtain
\begin{align}\label{eq:partialpN}
\begin{aligned}
\partial p_N(A) 
&=\frac{1}{M_V}\sum_{k=0}^{N} \partial(R_A^k)= -\frac{1}{M_V^2} \sum_{k=1}^{N} \sum_{j=0}^{k-1} R_A^{\,j} \otimes R_A^{\,k-1-j}\\
&= -\frac{1}{M_V^2} \sum_{\substack{j,m \geq 0 \\ j+m \leq N-1}} R_A^{\,j} \otimes R_A^{\,m}.
\end{aligned}
\end{align}
Thus, it follows from \eqref{eq:sharp} that for any \(H \in \mathcal{A}\),
\[
\partial p_N(A)\sharp H = -\frac{1}{M_V^2} \sum_{\substack{j,m \geq 0 \\ j+m \leq N-1}} R_A^{\,j} \, H \, R_A^{\,m}.
\]
Moreover, since \(A^{-1} = \frac{1}{M_V}\sum_{j=0}^{\infty} R_{A}^{\,j}\), we have
\[
-A^{-1} H A^{-1} = -\frac{1}{M_V^2} \sum_{j,m=0}^{\infty} R_A^{\,j} \, H \, R_{A}^{\,m}.
\]
According to \(\|R_A^{\,j} H R_A^{\,m}\| \leq \rho_T^{j+m}\|H\|\) and counting \((n+1)\) pairs \((j,m)\) with \(j+m = n\), we obtain
\begin{equation}\label{eq:firstderivtail}
\begin{aligned}
\bigl\|\partial p_N(A)\sharp H + A^{-1} H A^{-1}\bigr\|\leq \frac{\|H\|}{M_V^2} \sum_{n=N}^{\infty} (n+1)\rho_T^{n}.
\end{aligned}
\end{equation}
Then, taking \(A = V_s\) and \(H = f(V_s)\) in \eqref{eq:firstderivtail} yields
\begin{equation}\label{eq:firstderivlim}
\sup_{0 \leq s \leq T} \bigl\|\partial p_N(V_s)\sharp f(V_s) + V_s^{-1} f(V_s) V_s^{-1}\bigr\| \longrightarrow 0 \quad \text{as } N \to \infty,
\end{equation}
where $-V_s^{-1} f(V_s) V_s^{-1}=-\frac{\xi}{2}(V_t^{-1})^3  + \frac{b}{2}V_t^{-1}.$

On the other hand, the second derivative of $R_A^k$ is given by 
\begin{align}
    \begin{aligned}
        \partial^2(R_A^k)&=2 (\partial \otimes \mathrm{id}) \partial R_A^k 
        =-\frac{2}{M_V}\sum_{j=0}^{k-1} (\partial R_A^{j}) \otimes R_A^{k-1-j}\\
        &=-\frac{2}{M_V}\sum_{j=0}^{k-1} \bigl(-\frac{1}{M_V}  \sum_{i=0}^{j-1} R_A^{i} \otimes R_A^{j-1-i}\bigl) \otimes R_A^{k-1-j}\\
        &=\frac{2}{M_V^2}\sum_{\substack{i,j\geq0\\i+j\leq k-2}}R_A^i\otimes R_A^j\otimes R_A^{k-2-i-j}.
    \nonumber
    \end{aligned}
\end{align}
Then, it follows from \eqref{eq:deltau1} that 
\begin{align}
    \begin{aligned}
        \Delta_{\frac{\sigma}{2}(\mathbf{1}\otimes\mathbf{1})}p_N(A)=\frac{\sigma^2}{2M_V^3}\sum_{k=2}^{N}\sum_{\substack{i,j,l \geq 0 \\ i+j+l=k-2}}R_A^i\varphi(R_A^j)R_A^{l}.
        \nonumber
    \end{aligned}
\end{align}
Using \(A^{-1} = \frac{1}{M_V}\sum_{i=0}^{\infty}R_A^{i}\) together with \(\varphi(A^{-1}) = \frac{1}{M_V}\sum_{j=0}^{\infty}\varphi(R_A^{j})\) and sending \(N \to \infty\) gives the limit
\[
\frac{\sigma^2}{2M_V^3} \sum_{i,j,l=0}^{\infty} R_A^{\,i} \, \varphi(R_A^{\,j}) \, R_A^{\,l}
= \frac{\sigma^2}{2}A^{-1}\,\varphi(A^{-1})\,A^{-1}
= \frac{\sigma^2}{2}\varphi(A^{-1})A^{-2}.
\]
The truncation error involves all triples with \(i+j+l \geq N-1\). Using \(|\varphi(R_A^{\,j})| \leq \|R_A^{\,j}\| \leq \rho_T^{\,j}\) and counting \(\frac{(n+1)(n+2)}{2}\) triples \((i, j, l)\) with \(i+j+l = n\),
\begin{equation}
\begin{aligned}
\bigl\|\Delta_{\frac{\sigma}{2}(\mathbf{1}\otimes\mathbf{1})}p_N(A) - \frac{\sigma^2}{2}\varphi(A^{-1})A^{-2}\bigr\|
&\leq \frac{\sigma^2}{2M_V^3} \sum_{n=N-1}^{\infty} \frac{(n+1)(n+2)}{2}\,\rho_T^{n},
\nonumber
\end{aligned}
\end{equation}
which tends to $0$ uniformly on \(\{A : \operatorname{spec}(A) \subseteq[\ell_T, M_T]\}\). Hence
\begin{equation}\label{eq:secondderivlim}
\sup_{0 \leq s \leq T} \bigl\|\frac12\Delta_{\frac{\sigma}{2}(\mathbf{1} \otimes \mathbf{1})}p_N(V_s) - \frac{\sigma^2}{4}\varphi(V_s^{-1})V_s^{-2}\bigr\| \longrightarrow 0 \quad \text{as } N \to \infty.
\end{equation}

By \eqref{eq:vvv} and \eqref{eq:partialpN}, the free stochastic integrations in \eqref{eq:inverse} and \eqref{eq:ito-pN} are
\[
\frac{\sigma}{2}\,\partial p_N(V_s)\sharp \mathrm{d}W_s
= -\frac{\sigma}{2M_V^2} \sum_{\substack{j,m \geq 0 \\ j+m \leq N-1}} R_s^{\,j}\,\mathrm{d}W_s\,R_s^{\,m}
\]
and
\[
-\frac{\sigma}{2} V_s^{-1}\,\mathrm{d}W_s\,V_s^{-1}
= -\frac{\sigma}{2M_V^2} \sum_{j,m=0}^{\infty} R_s^{\,j}\,\mathrm{d}W_s\,R_s^{\,m},
\]
respectively. So we have
\[
\frac{\sigma}{2}\int_0^t \partial p_N(V_s)\sharp \mathrm{d}W_s + \frac{\sigma}{2}\int_0^t V_s^{-1}\,\mathrm{d}W_s\,V_s^{-1}
= -\frac{\sigma}{2M_V^2} \sum_{\substack{j,m \geq 0 \\ j+m \geq N}} \int_0^t R_s^{\,j}\,\mathrm{d}W_s\,R_s^{\,m}.
\]
Applying the free BG inequality, we obtain
\[
\Bigl\|\int_0^t R_s^{\,j}\,\mathrm{d}W_s\,R_s^{\,m}\Bigr\|
\leq 2\sqrt{2}\Bigl(\int_0^t \|R_s^{\,j}\|^2 \|R_s^{\,m}\|^2 \mathrm{d}s\Bigr)^{\frac12}
\leq 2\sqrt{2T}\,\rho_T^{\,j+m}.
\]
Therefore,
\begin{align}
\begin{aligned}
&\Bigl\|\frac{\sigma}{2}\int_0^t \partial p_N(V_s)\sharp \mathrm{d}W_s + \frac{\sigma}{2}\int_0^t V_s^{-1}\,\mathrm{d}W_s\,V_s^{-1}\Bigr\|\\
&\leq \frac{\sigma\sqrt{2T}}{M_V^2} \sum_{\substack{j,m \geq 0,\\  j+m \geq N}} \rho_T^{\,j+m}
= \frac{\sigma\sqrt{2T}}{M_V^2} \sum_{n=N}^{\infty} (n+1)\rho_T^{n}, 
\nonumber
\end{aligned}
\end{align}
which tends to $0$ as \(N \to \infty\). Hence
\begin{equation}\label{eq:stochconv}
\frac{\sigma}{2} \int_0^t \partial p_N(V_s)\sharp \mathrm{d}W_s \longrightarrow -\frac{\sigma}{2} \int_0^t V_s^{-1}\,\mathrm{d}W_s\,V_s^{-1}
\end{equation}
in operator norm, uniformly in  $t \in [0,T]$.

Define the limiting drift
\[
F(s) := -V_s^{-1} f(V_s) V_s^{-1} + \frac{\sigma^2}{4}\varphi(V_s^{-1})(V_s^{-1})^2,
\]
and its finite-\(N\) approximant \(F_N(s) := \partial p_N(V_s)\sharp f(V_s) + \frac{1}{2}\Delta_{\frac{\sigma}{2}(1\otimes1)} p_N(V_s)\). From \eqref{eq:firstderivlim} and \eqref{eq:secondderivlim} we have \(\sup_{0 \leq s \leq T}\|F_N(s) - F(s)\| \to 0\), and consequently
\[
\Bigl\|\int_0^t \bigl(F_N(s) - F(s)\bigr)\mathrm{d}s\Bigr\|
\leq \int_0^t \|F_N(s) - F(s)\|\mathrm{d}s
\leq T\sup_{0 \leq s \leq T}\|F_N(s) - F(s)\| \longrightarrow 0.
\]
Together with \eqref{eq:polyconv} and \eqref{eq:stochconv}, we may let \(N \to \infty\) in \eqref{eq:ito-pN} to obtain the limiting integral equation in operator norm:
\begin{equation}\label{eq:inverselimit}
V_t^{-1} = V_0^{-1} + \int_0^t \Bigl(-V_s^{-1} f(V_s) V_s^{-1} + \frac{\sigma^2}{4} \varphi(V_s^{-1})(V_s^{-1})^2\Bigr) \mathrm{d}s - \frac{\sigma}{2} \int_0^t V_s^{-1}\,\mathrm{d}W_s\,V_s^{-1}.
\nonumber
\end{equation}
Thus \(\{V_t^{-1}\}_{0 \leq t \leq T}\) is a free It\^o process and satisfies the following differential form
\[
\mathrm{d}(V_t^{-1}) = \Bigl(-\frac{\xi}{2}(V_t^{-1})^3 + \frac{\sigma^2}{4}\varphi(V_t^{-1})(V_t^{-1})^2 + \frac{b}{2}V_t^{-1}\Bigr)\mathrm{d}t - \frac{\sigma}{2}V_t^{-1}\,\mathrm{d}W_t\,V_t^{-1},
\]
which means \(\{V_t^{-1}\}_{0 \leq t \leq T}\) is a positive  solution of \eqref{eq:inverse}. Since \(T > 0\) was arbitrary, the result holds on every finite time interval.
\end{proof}

\section{Free analogue of the drift-implicit square-root Euler method}
\label{sec:method}
In this section, we introduce a drift-implicit square-root Euler scheme specifically for the free CIR equation and prove that it preserves positivity.

\begin{definition}\label{deffi:numer}
Fix a terminal time $T > 0$ and an integer $L \geq 1$. Define the uniform time grid $t_k = k \Delta t$ with step size $\Delta t = \frac{T}{L}$ for $k = 0, \dots, L$, the free drift-implicit square-root Euler scheme $\{\bar{V}_{t_{k}}\}_{k=0}^{L}$ is given by:
\begin{align}\label{eq:freeCIRfem}
\begin{aligned}
\bar{V}_{t_{k+1}}=\bar{V}_{t_k}+f(\bar{V}_{t_{k+1}})\Delta t+\frac{\sigma }{2}\Delta W_{t_{k}},\quad k=0,\ldots,L-1,
\end{aligned}
\end{align}
where $\Delta W_{t_{k}}:=W_{t_{k+1}}-W_{t_k}$, with initial value $ \bar{V}_0=V_0=\sqrt{X_0} \in \mathcal{A}_{+}$. The corresponding approximation to the free CIR equation \eqref{eq:freeCIRmain} is then given by \(\bar{X}_{t_{k}}=(\bar {V}_{t_{k}})^2.\) 
\end{definition}

The following theorem shows that \Cref{deffi:numer} is well-defined and establishes the existence and uniqueness of the numerical solution that preserves the positivity.

\begin{theorem}\label{prop:wellposed}
Given $\bar{V}_{0}\in\mathcal{A}_+$, then the scheme \eqref{eq:freeCIRfem} admits a unique solution
  $\{\bar{V}_{t_{k}}\}_{k=0}^{L}$ in $\mathcal{A}_+$. Moreover,  $(\bar{V}_{t_{k}})^{-1}$ is well-defined in $\mathcal{A}_+$, and the discrete solution
  $\bar{V}_{t_k}$ and $\bar{X}_{t_k}$ is adapted to $\mathcal{A}_{t_k}$,  $ k=0,\cdots,L$.
\end{theorem}
\begin{proof}
 For any $k=0,\cdots,L-1$ and $\bar{V}_{t_k}\in\mathcal{A}_+$,  set
  \begin{equation}\label{eq:notation-V-k-0}
  \zeta:=1+\frac{b}{2}\Delta t, \quad B:=\bar{V}_{t_k}+\tfrac{\sigma}{2}\Delta
  W_{t_{k}},\quad c:=\frac{\xi}{2}\Delta t.
  \end{equation}
Since $B$ is self-adjoint and $\zeta,c>0$, we have $B^{2}+4\zeta c\,\mathbf{1}>0.$
Thus, the positive square root
\[
S:=\sqrt{B^{2}+4\zeta c\,\mathbf{1}}
\]
is well-defined by the continuous functional calculus and commutes with $B$. Define
\begin{equation}\label{eq:notation-V-k-1}
\bar{V}_{t_{k+1}}:=\frac{B+S}{2\zeta}.
\end{equation}
We will show that the above $\bar{V}_{t_{k+1}}$ is strictly positive and thus invertible, and $\bar{V}_{t_{k+1}}$ satisfies  \eqref{eq:freeCIRfem}.

By a direct computation, we have 
\begin{align*}
\zeta \bar{V}_{t_{k+1}}^{2}-B\bar{V}_{t_{k+1}}
&=\frac{(B+S)^{2}}{4\zeta}-\frac{B(B+S)}{2\zeta}=\frac{(B+S)^{2}-2B(B+S)}{4\zeta}\\
&=\frac{S^{2}-B^{2}}{4\zeta}
=\frac{4\zeta c\,\mathbf{1}}{4\zeta }=c\mathbf{1}.
\end{align*}
Note that
\[
(B+S)(S-B)=S^{2}-B^{2}=4\zeta c\,\mathbf{1}.
\]
Since
$(B+S)^{-1}=\frac{S-B}{4\zeta c}.$
Hence $\bar{V}_{t_{k+1}}$ is invertible and
$
\bar{V}_{t_{k+1}}^{-1}=\frac{S-B}{2c}.
$
For every $\lambda\in\operatorname{spec}(B)$,
\[
f(\lambda)=\frac{\lambda+\sqrt{\lambda^{2}+4\zeta c}}{2\zeta}>0,
\]
which implies that $\bar{V}_{t_{k+1}}^{-1}$ is strictly positive. Consequently,  $(\bar{V}_{t_{k+1}})^{-1}$ is well-defined in $\mathcal{A}_+$ and  $\zeta\bar{V}_{t_{k+1}}^{2}-B\bar{V}_{t_{k+1}}-c\mathbf{1}=0$ implies \eqref{eq:freeCIRfem}.
 
For the uniqueness, let $W\in \mathcal{A}_+$ be another positive solution of \eqref{eq:freeCIRfem} produced by $\bar{V}_{t_{k}}$. Through straightforward calculation,  we have $[W,B]:=WB-BW=0$. The unital $C^{*}$‑algebra $C^{*}(\mathbf{1},B,W)$ is commutative. By the Gelfand representation, $B$ and $W$ correspond to real-valued continuous functions $\beta(\omega)$ and $w(\omega)$ on the spectrum of $C^{*}(\mathbf{1},B,W)$, respectively. Moreover, $w(\omega)\geq 0$, so the equation reduces pointwise to
\[
\zeta\,w(\omega)^{2}-\beta(\omega)\,w(\omega)-c=0.
\]
Since $\zeta,c>0$, the scalar roots $w_{+}(\omega)>0$ and $w_{-}(\omega)<0$, the unique non‑negative root reads
\[
w(\omega)=\frac{\beta(\omega)+\sqrt{\beta(\omega)^{2}+4\zeta c}}{2\zeta},
\]
so
$
W=\frac{B+\sqrt{B^{2}+4\zeta c\,\mathbf{1}}}{2\zeta}=V.
$
Thus the positive solution is unique.

Finally, we focus on the adaptedness.
Assume
$
\bar{V}_{t_{k}}\in\mathcal{A}_{t_{k}}\subseteq\mathcal{A}_{t_{k+1}}, \Delta W_{t_{k}}\in\mathcal{A}_{t_{k+1}}.
$
Then $B\in\mathcal{A}_{t_{k+1}}$. By the continuous functional calculus,
\[
\sqrt{B^{2}+4\zeta c\,\mathbf{1}}\in\mathcal{A}_{t_{k+1}},
\]
which yields $\bar{V}_{t_{k+1}}\in\mathcal{A}_{t_{k+1}}$. Hence the discrete solution $\bar{X}_{t_{k+1}}=(\bar{V}_{t_{k+1}})^2$ is adapted to $\mathcal{A}_{t_{k+1}}$ for any $k=0,\cdots,L-1$.
\end{proof}

\subsection{The operator norm bounds for the numerical solution}
In this section, we prove  two operator norm bounds for the numerical solution. The first bound holds for sufficiently small step sizes on a fixed interval, which enables us to establish the strong convergence of the proposed scheme. The second bound holds for every fixed $\Delta t>0$ (without any step size restriction) on all time steps $k$, which facilitates the derivation of the long-time  exponential stability for arbitrary step sizes. 

\begin{lemma}
\label{lem:numerical_bound}
Fix \(T>0\) and let \(t_k=k\Delta t\), \(k=0,\ldots,L\) with \(T=L\Delta t\). Let \(V_t\) and \(\bar V_{t_k}\) denote the exact and numerical solutions of \eqref{eq:freeCIRlam} and \eqref{eq:freeCIRfem} respectively, with $V_0=\bar V_{t_0}$. Then there exist constants \(h_0>0\) and \(K'_T>0\), independent of \(k\), such that, whenever \(\Delta t\leq h_0\),
\[
\max_{0\leq k\leq L}
\|V_{t_k}-\bar V_{t_k}\|
\leq K'_T(\Delta t)^{\frac12}.
\]
Moreover, we have
\[\frac{\ell_T}{2}\mathbf{1}\leq\bar V_{t_k}\leq\left(M_V+\frac{\ell_T}{2}\right)\mathbf{1},\]
where \(\ell_T\) and \(M_V\) are given in \Cref{lem:V_bound}.
\end{lemma}

\begin{proof}
By \Cref{lem:V_bound}, for \(0\leq s,t\leq T\), we have
\[
\ell_T\mathbf{1}\leq V_t\leq M_V\mathbf{1},
\quad
\|V_t-V_s\|\leq K_T|t-s|^{\frac12}.
\]
Moreover, for \(A,B\geq \ell_T\mathbf{1}\),
\[
\|f(A)-f(B)\|
\leq
L_0\|A-B\|,
\quad
L_0:=\frac{\xi}{2\ell_T^2}+\frac{b}{2},
\]
where we have used $A^{-1}-B^{-1}=-A^{-1}(A-B)B^{-1}.$ For later use, define the remainder term
\begin{align}\label{eq:rk}
    \begin{aligned}
        r_k:=\int_{t_k}^{t_{k+1}}\bigl(f(V_s)-f(V_{t_{k+1}})\bigr)\,\mathrm{d}s, \quad k=0,\ldots,L-1.
    \end{aligned}
\end{align}
Then
\begin{align}
  \begin{aligned}
\|r_k\|
&\leq\int_{t_k}^{t_{k+1}}
\|f(V_s)-f(V_{t_{k+1}})\|\,\mathrm{d}s\\
&\leq
L_0K_T
\int_{t_k}^{t_{k+1}}
|t_{k+1}-s|^{\frac12}\,\mathrm{d}s =
\frac{2}{3}L_0K_T(\Delta t)^{\frac32}.
\nonumber
\end{aligned}  
\end{align}
Set $K_r:=\frac{2}{3}L_0K_T$ and $L_1
:=
\frac{2\xi}{\ell_T^2}
+\frac{b}{2}$, and choose \(h_0>0\)
sufficiently small such that
\[
h_0L_1\le \frac12\quad  \text{and}
\quad
K_rT e^{2L_1T}h_0^{\frac12}\le \frac{\ell_T}{4}.
\]
Then, for every $0<\Delta t\le h_0$, define
\[
\rho:=(1-\Delta t L_1)^{-1},
\quad
E_0:=0,
\quad
E_{k+1}
:=
\rho\bigl(E_k+K_r(\Delta t)^{\frac32}\bigr).
\]
Since $\rho\leq 1+2\Delta t L_1\leq e^{2L_1\Delta t},$ we obtain, for \(k\leq L,\) 
\begin{align}
\begin{aligned}
E_k
&=
K_r(\Delta t)^{\frac32}
\sum_{j=1}^{k}\rho^j \leq K_r(\Delta t)^{\frac32}k\rho^k\leq
K_rT e^{2L_1T}(\Delta t)^{\frac12}\leq\frac{\ell_T}{4}.
\nonumber
\end{aligned}
\end{align}

We now prove by induction that $\|V_{t_k}-\bar V_{t_k}\|\leq E_k.$ The assertion is immediate for \(k=0\). Suppose that it holds for some
\(k<L\), consider the closed ball
$
\mathcal{B}
:=
\left\{
A\in\mathcal{A}_{sa}:
\|A-V_{t_{k+1}}\|\leq\frac{\ell_T}{2}
\right\}.
$
Since \(V_{t_{k+1}}\geq\ell_T\mathbf{1}\), any \(A\in\mathcal{B}\) satisfies
$A\geq\frac{\ell_T}{2}\mathbf{1},$
which implies \(\mathcal B \subset \mathcal A_+\). Moreover, for \(A,B\in\mathcal{B}\), we have
\[
\|f(A)-f(B)\|
\leq
L_1\|A-B\|.
\]
Define
\[
\Phi(A):=\bar V_{t_k}+\frac{\sigma}{2}\Delta W_{t_{k}}+\Delta t\,f(A).
\]
Using the integral form of \eqref{eq:freeCIRlam} and the definition of $r_k$,  for \(A\in\mathcal{B}\), we obtain
\begin{align}
    \begin{aligned}
    \|\Phi(A)-V_{t_{k+1}}\|&=\bigg\|\bigl(\bar V_{t_k}+\frac{\sigma}{2}\Delta W_{t_{k}}+\Delta t\,f(A)\bigl)-\bigl(V_{t_k}+\Delta t\,f(V_{s})+\int_{t_k}^{t_{k+1}}\mathrm{d}W_s\bigl)\bigg\|\\
    &=\|\bar V_{t_k}+\Delta t\,f(A)-V_{t_k}-\Delta t\,f(V_{t_{k+1}})-r_k\|\\
    &\leq
    \|V_{t_k}-\bar V_{t_k}\|
    +\|r_k\|
    +\Delta t\|f(A)-f(V_{t_{k+1}})\| \\
    &\leq
    E_k+K_r(\Delta t)^{\frac32}
    +\Delta t L_1\|A-V_{t_{k+1}}\| \\
    &\leq
    (1-\Delta t L_1)E_{k+1}
    +\frac{\Delta t L_1\ell_T}{2} \\
    &\leq
    \frac{(1-\Delta t L_1)\ell_T}{4}
    +\frac{\Delta t L_1\ell_T}{2}
    \leq
    \frac{\ell_T}{2}.
    \nonumber
    \end{aligned}
\end{align}
Thus \(\Phi(\mathcal{B})\subseteq \mathcal{B}\). 
Since \(\Delta tL_1<1\) and
\[
\|\Phi(A)-\Phi(B)\|\leq\Delta t \|f(A)-f(B)\|\leq  L_1\Delta t\|A-B\|,
\]
then the map \(\Phi\) is a contraction on \(\mathcal B\).
Therefore, by the uniqueness of the positive solution established in
\Cref{prop:wellposed} and the Banach fixed-point theorem, we get
\[
\bar V_{t_{k+1}}=\Phi(\bar V_{t_{k+1}})\in \mathcal{B}.
\]
Consequently, 
\[
\|\bar V_{t_{k+1}}-V_{t_{k+1}}\|
\leq
E_k+K_r(\Delta t)^{\frac32}
+\Delta tL_1
\|\bar V_{t_{k+1}}-V_{t_{k+1}}\|.
\]
Moving the last term on the right-hand side to the left-hand side, we have
\[
\|\bar V_{t_{k+1}}-V_{t_{k+1}}\|
\leq \frac{E_k+K_r(\Delta t)^{\frac32}}{1-\Delta tL_1}=
E_{k+1}.
\]
Then
$\max_{0\leq k\leq L}
\|V_{t_k}-\bar V_{t_k}\|
\leq
K'_T(\Delta t)^{\frac12},\ K'_T:=K_rT e^{2L_1T}.$

Finally, decreasing \(h_0\) once more if necessary, we may assume that
$K'_T(\Delta t)^{\frac12}
\leq
\frac{\ell_T}{2}.$
Since $\ell_T\mathbf{1}\leq V_{t_k} \leq M_V\mathbf{1},$ it follows that
\[
\frac{\ell_T}{2}\mathbf{1}
\leq
\bar V_{t_k}
\leq
\left(M_V+\frac{\ell_T}{2}\right)\mathbf{1},
\quad
0\leq k\leq L.
\]
\end{proof}

\begin{lemma}\label{lem:num_global_bound}
Let $\{\bar V_{t_k}\}_{k\ge0}$ be the positive solution of
\eqref{eq:freeCIRfem}. Then we have
\[
\sup_{k\ge0}\|\bar V_{t_k}\|
\le \frac{\sigma}{\sqrt b}
+\sqrt{
\|\bar X_0\|
+\frac{a+7\sigma^2/4}{b}
}.
\]
\end{lemma}
\begin{proof}
Rewriting the numerical scheme gives
\[
\bar V_{t_k}+\frac{\sigma}{2}\Delta W_{t_k}
=
\bar V_{t_{k+1}}
-\Delta t\,f(\bar V_{t_{k+1}}).
\]
Since
\(
Vf(V)=f(V)V
=\frac{\xi}{2}\mathbf{1}-\frac{b}{2}V^2,
\)
squaring both sides, we have
\begin{align*}
&\bar X_{t_k}+\frac{\sigma^2}{4}(\Delta W_{t_k})^2+\frac{\sigma}{2}
\bigl(
\bar V_{t_k}\Delta W_{t_k}
+\Delta W_{t_k}\bar V_{t_k}
\bigr)\\
&=(1+b\Delta t)\bar X_{t_{k+1}}-\xi\Delta t\,\mathbf{1}+(\Delta t)^2 f(\bar V_{t_{k+1}})^2.
\end{align*}
Rearranging yield
\begin{align*}
(1+b\Delta t)\bar X_{t_{k+1}}
={}&\bar X_{t_k}+\xi\Delta t\,\mathbf{1}
+\frac{\sigma}{2}
\bigl(
\bar V_{t_k}\Delta W_{t_k}
+\Delta W_{t_k}\bar V_{t_k}
\bigr)\\
&+\frac{\sigma^2}{4}(\Delta W_{t_k})^2
-(\Delta t)^2 f(\bar V_{t_{k+1}})^2.
\end{align*}
Because $f(\bar V_{t_{k+1}})$ is self-adjoint,
its square is nonnegative. Moreover,
$\|\Delta W_{t_k}\|=2\sqrt{\Delta t}$ implies
\[
(\Delta W_{t_k})^2\le 4\Delta t\,\mathbf{1}.
\]
Set $\gamma=(1+b\Delta t)^{-1}\in(0,1)$.
Strict positivity of the numerical solution and iteration therefore gives, 
\begin{equation}\label{eq:numerical-operator-upper}
\begin{aligned}
0\le\bar X_{t_k}
&\le
\gamma^k\bar X_0
+(\xi+\sigma^2)\Delta t
\sum_{i=1}^{k}\gamma^i\,\mathbf{1}
+J_k\\
&=
\gamma^k\bar X_0
+\frac{\xi+\sigma^2}{b}
(1-\gamma^k)\mathbf{1}
+J_k,
\end{aligned}
\end{equation}
where
\(
J_k:=
\frac{\sigma}{2}
\sum_{j=0}^{k-1}\gamma^{k-j}
\bigl(
\bar V_{t_j}\Delta W_{t_j}
+\Delta W_{t_j}\bar V_{t_j}
\bigr).
\)
Fix an arbitrary integer $K\ge1$ and define
\[
S_K:=\max_{0\le j\le K}\|\bar X_{t_j}\|<\infty.
\]
For each $1\le k\le K$, introduce the step process
\[
H_s:=\gamma^{k-j}\bar V_{t_j},
\qquad
s\in[t_j,t_{j+1}),\quad j=0,\ldots,k-1.
\]
Since
$\bar V_{t_j}\in\mathcal A_{t_j}$, and
\(
J_k=
\frac{\sigma}{2}
\left(
\int_0^{t_k}H_s\,dW_s
+\int_0^{t_k}dW_s\,H_s
\right).
\)
Applying the free BG inequality
in the operator norm to these two integrals gives
\begin{align*}
\|J_k\|
&\le
2\sqrt{2}\,\sigma
\left(
\int_0^{t_k}\|H_s\|^2\,ds
\right)^{\frac12}=
2\sqrt{2}\,\sigma
\left(
\Delta t
\sum_{j=0}^{k-1}
\gamma^{2(k-j)}\|\bar V_{t_j}\|^2
\right)^{\frac12}\\
&\le
2\sqrt{2}\,\sigma\sqrt{S_K}
\left(
\Delta t\sum_{m=1}^{k}\gamma^{2m}
\right)^{\frac12},
\end{align*}
where the last inequality is due to
$\|\bar V_{t_j}\|^2=\|\bar X_{t_j}\|\le S_K$.
The geometric sum satisfies
\[
\Delta t\sum_{m=1}^{k}\gamma^{2m}
\le
\frac{\Delta t\,\gamma^2}{1-\gamma^2}
=
\frac{1}{b(2+b\Delta t)}
\le\frac{1}{2b}.
\]
Consequently,
\[
\|J_k\|
\le\frac{2\sigma}{\sqrt b}\sqrt{S_K},
\qquad 1\le k\le K.
\]

According to
\eqref{eq:numerical-operator-upper}, we obtain
\[
\|\bar X_{t_k}\|
\le
\|\bar X_0\|
+\frac{\xi+\sigma^2}{b}
+\frac{2\sigma}{\sqrt b}\sqrt{S_K},
\qquad 1\le k\le K.
\]
Taking the maximum over $0\le k\le K$ yields
\[
S_K
\le
\|\bar X_0\|
+\frac{\xi+\sigma^2}{b}
+\frac{2\sigma}{\sqrt b}\sqrt{S_K}.
\]
Solving this quadratic inequality in $\sqrt{S_K}$, we conclude that
\[
\|\bar V_{t_k}\|
\le\sqrt{S_K}
\le
\frac{\sigma}{\sqrt b}
+\sqrt{
\|\bar X_0\|
+\frac{a+7\sigma^2/4}{b}
},
\qquad 0\le k\le K.
\]
The right-hand side is independent of $\Delta t$ and $K$.
Since $K$ is arbitrary, the bound holds for every
$\Delta t>0$ and every $k\ge0$.
\end{proof}

\section{Convergence and stability of the  numerical method} 
\label{sec:results}
In this section, we prove the main results about convergence and stability  for the free drift-implicit square-root Euler method.

\subsection{Strong convergence}\label{s:5.1}
We establish the strong convergence of the free drift-implicit square-root Euler method in \Cref{thm:convergence}.
We emphasize that our error analysis is based on a decomposition of the remainder term \(r_k\) that isolates a martingale component, which allows us to improve the classical results.
\begin{theorem}
\label{thm:convergence}
Fix \(T>0\), define the grid points \(t_k=k\Delta t\) for \(k=0,\ldots,L\), where
\(T=L\Delta t\). Let \(V_t\) be the exact solution of \eqref{eq:freeCIRlam} and \(\bar V_{t_k}\) be the numerical approximation generated by scheme of \eqref{eq:freeCIRfem} with $ V_0=\bar V_0.$ 
Let $ X_t$ be the positive solution of the free CIR equation \eqref{eq:freeCIRmain} satisfying \eqref{eq:square root}, and  let $\bar X_{t_k}$ be defined as in  \Cref{deffi:numer}. Then there exist constants \(C_T>0\) and \(h_0>0\), independent of
\(\Delta t\) and \(L\),  such that for every \(0<\Delta t\le h_0\),
\[
\max_{0\leq k\leq L}
\|X_{t_k}-\bar X_{t_k}\|_2
\leq
C_T\Delta t.
\]
\end{theorem}

To prove \Cref{thm:convergence}, we need the following lemmas.
\begin{lemma}
\label{lem:localerror}
For any $ k=0,\ldots,L-1$ and the remainder term \(r_k\) defined in the proof of
\Cref{lem:numerical_bound}, there exist self-adjoint elements  $A_k,M_k\in\mathcal{A}_{t_{k+1}}$  
such that
\[
r_k=A_k+M_k.
\]
Moreover, there exist constants \(C_A,C_M>0\), independent of
\(k\), \(\Delta t\) and \(L\), such that
\[
\|A_k\|_2
\leq
C_A(\Delta t)^2,
\quad
\|M_k\|_2^2
\leq
C_M(\Delta t)^3, \quad \mathbb{E}_{t_k}(M_k)=0,
\]
where
\(\mathbb{E}_{t_k}:\mathcal{A}\rightarrow\mathcal{A}_{t_k}\)
denotes the trace-preserving conditional expectation.
\end{lemma}
\begin{proof}
Recall that $f(V_t)=\frac{\xi}{2}V_t^{-1}-\frac{b}{2}V_t.$ \Cref{prop:inverse} and the free SDE \eqref{eq:freeCIRlam} give for any $t\in[0,T],$
\begin{align}
    \begin{aligned}
        \mathrm{d}f(V_t)=G_t\mathrm{d}t-\frac{\xi\sigma }{4}V_t^{-1}\mathrm{d}W_tV_t^{-1}-\frac{b\sigma}{4}\mathrm{d}W_t,
        \nonumber
    \end{aligned}
\end{align}
where
\(
G_t=-\frac{\xi^2}{4}(V_t^{-1})^3+
\frac{\xi\sigma^2}{8}\varphi(V_t^{-1})(V_t^{-1})^2+
\frac{b^2}{4}V_t.
\)
For \(s\in[t_k,t_{k+1}]\), integrating the above equation from \(s\) to \(t_{k+1}\) and then integrating with respect to \(s\) over \([t_k,t_{k+1}]\) give $r_k=A_k+M_k,$  where
\[
A_k:=-\int_{t_k}^{t_{k+1}}
(u-t_k)G_u\mathrm{d}u,
\]
and
\[
M_k:=\frac{\xi\sigma}{4}
\int_{t_k}^{t_{k+1}}
(u-t_k)
V_u^{-1}\mathrm{d}W_uV_u^{-1} +\frac{b\sigma}{4}
\int_{t_k}^{t_{k+1}}
(u-t_k)\mathrm{d}W_u.
\]
Obviously, \(A_k\) and \(M_k\) are self-adjoint and belong to
\(\mathcal{A}_{t_{k+1}}\).

By \Cref{lem:V_bound}, there exist
constants \(G_*,C_*>0\) independent of \(k\) and \(\Delta t\), such that
\[
\sup_{0\leq t\leq T}\|G_t\|\leq G_*:=\frac{2\xi^2+\xi\sigma}{8\ell_T^3}+\frac{b^2}{4}M_V,
\]
and 
\[
\frac{\xi\sigma}{4}\sup_{0\leq t\leq T}\|V_t^{-1}\|^2+\frac{b\sigma}{4}\leq C_*:=\frac{\xi\sigma}{4\ell_T^2}+\frac{b\sigma}{4}.
\]
For the drift part $A_k$, it follows from $L^\infty(\mathcal{A}) \subset L^2(\mathcal{A})$ that
\[
\|G_u\|_2 \leq \|G_u\| \leq G_*,
\]
we obtain
\[
\|A_k\|_2\leq
\int_{t_k}^{t_{k+1}}
(u-t_k)\|G_u\|_2\mathrm{d}u\leq G_*\int_{t_k}^{t_{k+1}}
(u-t_k)\mathrm{d}u=
\frac{G_*}{2}(\Delta t)^2.
\]
Hence, with $C_A:=\frac{G_*}{2}$, we have
$\|A_k\|_2
\leq
C_A(\Delta t)^2.$
For the martingale part $M_k$, the free It\^o isometry 
yields
\[
\|M_k\|_2^2\leq C_*^2 \int_{t_k}^{t_{k+1}} (u-t_k)^2\mathrm{d}u =\frac{C_*^2}{3}(\Delta t)^3.
\]
Thus, with $C_M:=\frac{C_*^2}{3},$ we obtain
$\|M_k\|_2^2
\leq
C_M(\Delta t)^3.$

Finally, since \(M_k\) is a free stochastic integral over \([t_k,t_{k+1}]\), by the martingale property of the free stochastic integral, we get
\(
\mathbb{E}_{t_k}(M_k)=0.
\)
\end{proof}

\begin{lemma}\label{lem:dissipativity}
For any \(A,B\in\mathcal{A}_+\), we have
$$
\varphi\!\left((A-B)(f(A)-f(B))\right)
\leq
-\frac{b}{2}\|A-B\|_2^2.
$$
\end{lemma}

\begin{proof}
Denote $D:=A-B$. Using $A^{-1}-B^{-1}=-A^{-1}(A-B)B^{-1}$ and the tracial property of \(\varphi\), we obtain
\begin{align}
\begin{aligned}
\varphi\!\left(D(A^{-1}-B^{-1})\right)
=-\varphi\!\left(DA^{-1}DB^{-1}\right)=
-\|A^{-\frac12}DB^{-\frac12}\|_2^2
\leq 0.
\nonumber
\end{aligned}
\end{align}
Therefore, we have
\begin{align}
\begin{aligned}
\varphi\!\left(D(f(A)-f(B))\right)
=
\frac{\xi}{2}
\varphi\!\left(D(A^{-1}-B^{-1})\right)
-\frac{b}{2}\|D\|_2^2\leq
-\frac{b}{2}\|D\|_2^2.
\nonumber
\end{aligned}
\end{align}
\end{proof}

Now, We are in the position to prove \Cref{thm:convergence}.
\begin{proof}[Proof of \Cref{thm:convergence}]
From \eqref{eq:freeCIRlam} and \eqref{eq:freeCIRfem}, it follows that
\begin{align}\label{eq:pointerror}
\begin{aligned}
V_{t_{k+1}}-\bar{V}_{t_{k+1}}&=V_{t_{k}}-\bar{V}_{t_{k}}+\int_{t_k}^{t_{k+1}}\big(f(V_{s})-f(\bar V_{t_{k+1}})\big)\mathrm{d}s.
\end{aligned}
\end{align}
Define the error $e_{k} := V_{t_{k}} - \bar{V}_{t_{k}}$ and $\Delta e_k=e_{k+1}-e_k.$
Then \eqref{eq:pointerror} can be rewritten as
\begin{align}
\begin{aligned}
\Delta  e_{k}&=r_k+\big(f(V_{t_{k+1}})-f(\bar V_{t_{k+1}})\big)\Delta t=A_k+M_k+\big(f(V_{t_{k+1}})-f(\bar V_{t_{k+1}})\big)\Delta t,
\nonumber
\end{aligned}
\end{align}
where $A_k,M_k$ come from \Cref{lem:localerror}.
Now, consider the inner product term:
\begin{align}\label{eq:ek+1ek}
\begin{aligned}
\varphi (e_{k+1}\Delta  e_{k}) = \varphi (e_{k+1} A_k)+  \varphi (e_{k+1} M_k) + \varphi \bigg(e_{k+1}\big(f(V_{t_{k+1}})-f(\bar V_{t_{k+1}})\big)\Delta t\bigg).
\end{aligned}
\end{align}
On the one hand, since \(e_k\in\mathcal A_{t_k}\) and \(\mathbb{E}_{t_k}(M_k)=0\), the trace-preserving property of the conditional expectation gives
\[
\varphi(e_kM_k)=\varphi\bigl(e_k\mathbb{E}_{t_k}(M_k)\bigr)=0.
\]
Therefore, by the noncommutative Cauchy--Schwarz inequality,
\[
2\varphi(e_{k+1}M_k)=2\varphi((e_k+\Delta e_k)M_k)=2\varphi(\Delta e_kM_k)\leq\frac12\|\Delta e_k\|_2^2+2\|M_k\|_2^2;
\]
\[
2\varphi(e_{k+1}A_k)\leq\frac{b\Delta t}{2}\|e_{k+1}\|_2^2+\frac{2}{b\Delta t}\|A_k\|_2^2.
\]
By \Cref{lem:dissipativity}, we have
\[
2\varphi \bigg(e_{k+1}\big(f(V_{t_{k+1}})-f(\bar V_{t_{k+1}})\big)\Delta t\bigg)\leq -b\Delta t\|e_{k+1}\|_2^2.
\]
On the other hand, 
\[
2\varphi (e_{k+1}\Delta  e_{k})=\varphi(e_{k+1}^2-e_{k}^2+(\Delta  e_{k})^2)=\|e_{k+1}\|^2_2-\|e_{k}\|_2^2+\|\Delta  e_{k}\|_2^2.
\]
Combining the above terms, we obtain
\[
\|e_{k+1}\|^2_2-\|e_{k}\|_2^2+\|\Delta  e_{k}\|_2^2\leq -\frac{b\Delta t}{2}\|e_{k+1}\|_2^2+\frac{2}{b\Delta t}\|A_k\|_2^2+\frac12\|\Delta e_k\|_2^2+2\|M_k\|_2^2.
\]
Consequently,
\[
\left(
1+\frac{b\Delta t}{2}
\right)
\|e_{k+1}\|_2^2
+
\frac12\|\Delta e_k\|_2^2
\leq
\|e_k\|_2^2
+
\frac2{b\Delta t}\|A_k\|_2^2+2\|M_k\|_2^2\leq \|e_k\|_2^2+C (\Delta t)^3,
\]
where \(C>0\) is independent of \(\Delta t\). Dropping the nonnegative term $\frac12\|\Delta e_k\|_2^2$, setting $q:=(1+\frac{b\Delta t}{2})^{-1}<1$ and $e_{0}=0,$ we obtain, for $n=1,\ldots,L$,
\begin{align}
\begin{aligned}
\|e_n\|_2^2&\leq q\|e_{n-1}\|_2^2+Cq(\Delta t)^3\leq \dots \leq q^{n}\|e_{0}\|_2^2+Cq(\Delta t)^3\sum_{i=0}^{n-1}q^i\\
&\leq
Cq(\Delta t)^3\sum_{i=0}^{n-1}q^i\leq \frac{Cq}{1-q}(\Delta t)^3=\frac{2C}{b}(\Delta t)^2.
\nonumber
\end{aligned}    
\end{align}
Thus,
$
\max_{0\leq k\leq L}
\|V_{t_k}-\bar V_{t_k}\|_2
\leq\sqrt{{2C}/{b}}\Delta t.
$
The noncommutative H\"older inequality gives
\[
\|X_{t_k}-\bar X_{t_k}\|_2
\leq
\left(
\|V_{t_k}\|
+
\|\bar V_{t_k}\|
\right)
\|e_k\|_2.
\]
By \Cref{lem:V_bound,lem:numerical_bound}, the factor in parentheses is uniformly bounded for \(0\leq k\leq L\) and \(0<\Delta t\leq h_0\). Hence,
\[
\max_{0\leq k\leq L}
\|X_{t_k}-\bar X_{t_k}\|_2
\leq
C_T\Delta t.
\]
This completes the proof.
\end{proof}

\subsection{\texorpdfstring{Exponential stability in $L^2(\varphi)$}{Exponential stability in L2(phi)}}\label{s:4.3}
In this section, we investigate the long-time stability of the scheme \eqref{eq:freeCIRfem} with respect to perturbations of the initial value for both the exact and numerical solutions.
Analogue to the classical cases \cite{higham2003exponential, DAmbrosio2023,HUANG20124016}, we introduce the following definitions for both the continuous and discrete settings, respectively.
\begin{definition}\label{def:expstable}
Let $\{X_t\}_{t\ge 0}$ and $\{\tilde{X}_t\}_{t\ge 0}$ be two solutions of \eqref{eq:freeSDE} with initial values $X_0,\tilde{X}_0$. The free SDEs are said to be exponentially stable in $L^2(\varphi)$ if there exist constants $C,\lambda_1>0,$ independent of $t$, such that
\[
\|X_t-\tilde{X}_t\|^2_2 \le Ce^{-\lambda_1t}\|X_0-\tilde{X}_0\|^2_2 .
\]
\end{definition}

\begin{definition}\label{def:numstable}
For a given step size $\Delta t > 0$, let $\{\bar{X}_{t_k}\}_{k\in\mathbb{N}}$ and $\{\hat X_{t_k}\}_{k\in\mathbb{N}}$ be two numerical solutions to the free SDEs \eqref{eq:freeSDE} generated by a numerical method with initial values $\bar{X}_0, \hat X_0$, where $t_k = k\Delta t$.  The numerical method is said to be exponentially stable in $L^2(\varphi)$ if there exist constants $C,\lambda_2 > 0,$ independent of $k$ and  $\Delta t$, such that
\[
\|\bar{X}_{t_k} - \hat X_{t_k}\|^2_2 \le C e^{-\lambda_2k\Delta t}\|\bar{X}_{0} - \hat X_{0}\|^2_2 .
\]
\end{definition}

\begin{theorem}\label{thm:stability1}
  Let $X_t$ and $\tilde{X}_t$ be two solutions of the free CIR equation
  \eqref{eq:freeCIRmain} with initial values $X_0,\tilde{X}_0\in\mathcal{A}_+$.
Then 
\begin{equation}\label{eq:sta1}
\|X_t-\tilde{X}_t\|^2_2 \le Ce^{-bt}\|X_0-\tilde{X}_0\|^2_2 .
\nonumber
\end{equation}
\end{theorem}

  \begin{proof}
Let $V_t$ and $U_t$ be the solutions of \eqref{eq:freeCIRlam} with $V_0=\sqrt{X_0}$ and $U_0=\sqrt{\tilde{X}_0}$, i.e., $V_t^2=X_t$ and $U_t^2=\tilde{X}_t$. Subtracting the equations for $V_t$ and $U_t$, the additive noise $\tfrac{\sigma}{2}\,\mathrm{d}W_t$ cancels, so that
\[
  \mathrm{d}(V_t-U_t)
  =\bigl(f(V_t)-f(U_t)\bigr)\mathrm{d}t.
\]
Using the free It\^o product rule and applying $\varphi$ to the both sides, \Cref{lem:dissipativity} implies 
\begin{equation*}
    \begin{aligned}
         \frac{\mathrm{d}}{\mathrm{d}t}\,\varphi\bigl((V_t-U_t)^2\bigr)\leq 2\varphi((V_t-U_t)(f(V_t)-f(U_t))) \leq -b\,\varphi\bigl((V_t-U_t)^2\bigr),
    \end{aligned}
\end{equation*}
  so Gronwall's inequality yields
  \begin{equation}\label{eq:VUdecay}
  \|V_t-U_t\|_2^2\le e^{-bt}\|V_0-U_0\|_2^2.
  \end{equation}
By \Cref{lem:V_bound}, for $t\ge0$,
\[
\|V_t\|\le M_V,\quad
\|U_t\|\le M_U.
\]
Set
$M_*:=M_V+M_U.$
Using H\"older's inequality yields
\[
\|X_t-\widetilde X_t\|_2=\|V_t(V_t-U_t)+(V_t-U_t)U_t\|_2
\le M_*\|V_t-U_t\|_2.
\]
Moreover, set $m_*:=
\inf\operatorname{spec}(V_0)+\inf\operatorname{spec}(U_0)>0$ and denote $D_0:=V_0-U_0$ then
$V_0D_0+D_0U_0=X_0-\widetilde X_0.$
Since $V_0,U_0 \in \mathcal{A}_+$, so we have
\[
D_0
=
\int_0^\infty
e^{-sV_0}(X_0-\widetilde X_0)e^{-sU_0}\,\mathrm{d}s,
\]
and hence
\[
\|V_0-U_0\|_2
\le
\frac{1}{m_*}\|X_0-\widetilde X_0\|_2.
\]
Combining this with \eqref{eq:VUdecay}, we obtain
\(
\|X_t-\widetilde X_t\|^2_2
\le
\frac{M^2_*}{m^2_*}e^{-b t}
\|X_0-\widetilde X_0\|^2_2.
\)
  \end{proof}

\begin{theorem}\label{thm:numstability}
For a given step size $\Delta t > 0$, let $\{\bar{X}_{t_k}\}_{k\in\mathbb{N}}$ and $\{\hat X_{t_k}\}_{k\in\mathbb{N}}$ be the numerical solutions  to the free CIR equation \eqref{eq:freeCIRmain} produced by the free drift-implicit square-root Euler scheme \eqref{eq:freeCIRfem} with initial values $\bar{X}_0, \hat X_0 \in \mathcal{A}_+$.  
Then there exist constant $C> 0$ independent of $k$ and  $\Delta t$, such that
 \[
\|\bar{X}_{t_k} - \hat X_{t_k}\|^2_2 \le Ce^{-\frac{b}{1+b\Delta t}k\Delta t}\|\bar{X}_{0} - \hat X_{0}\|^2_2.
\]
\end{theorem}

\begin{proof}
Let $\bar{V}_{t_k}$ and $\bar{U}_{t_k}$ to be the square-root approximations corresponding to $\bar{X}_{t_k}$ and $\hat X_{t_k}$, respectively, and denote the difference $\bar{D}_{t_k} := \bar{V}_{t_k} - \bar{U}_{t_k}$. From the scheme definition \eqref{eq:freeCIRfem}, we have the recurrence
\begin{align}\label{eq:numVUrecur}
\bar{D}_{t_{k+1}} = \bar{D}_{t_k} + \big(f(\bar{V}_{t_{k+1}}) - f(\bar{U}_{t_{k+1}})\big) \Delta t.
\end{align}
Multiplying both sides of \eqref{eq:numVUrecur} on the left by $\bar{D}_{t_{k+1}}$, applying the trace $\varphi$ and \Cref{lem:dissipativity}, we obtain
\begin{align}\label{eq:numVUstep}
\begin{aligned}
\|\bar{D}_{t_{k+1}}\|_2^2
&= \varphi\bigl(\bar{D}_{t_{k+1}}\bar{D}_{t_k}\bigr)
   + \varphi \bigg(\bar{D}_{t_{k+1}}\big(f(\bar{V}_{t_{k+1}}) - f(\bar{U}_{t_{k+1}})\big)\bigg) \Delta t \\[2pt]
&\le \frac{1}{2}\|\bar{D}_{t_{k+1}}\|_2^2 + \frac{1}{2}\|\bar{D}_{t_k}\|_2^2 - \frac{b\Delta t}{2}\|\bar{D}_{t_{k+1}}\|_2^2.
\nonumber
\end{aligned}
\end{align}
Rearranging terms yields
\[
(1 + b\Delta t)\,\|\bar{D}_{t_{k+1}}\|_2^2 \le \|\bar{D}_{t_k}\|_2^2.
\]
Iterating this inequality $k$ times and using $1 - x \le e^{-x}$ for $x = b\Delta t/(1+b\Delta t) \in (0,1)$, we obtain the geometric decay
\begin{align}\label{eq:numVUdecay}
\|\bar{D}_{t_k}\|_2^2 \le \Bigl(\frac{1}{1+b\Delta t}\Bigr)^{\!k} \|\bar{D}_0\|_2^2
\le e^{-\frac{b\,k\Delta t}{1+b\Delta t}}\|\bar{D}_0\|_2^2.
\end{align}
Now consider the difference of the squared processes: $\bar{Z}_{t_k} := \bar{X}_{t_k} - \hat X_{t_k} =\bar V_{t_k}\bar D_{t_k}+\bar D_{t_k}\bar U_{t_k}.$
Using \Cref{lem:num_global_bound}, there exists a constant $\bar M_*,$ independent of $k$ and $\Delta t$, such that
\[
\|\bar V_{t_k}\|+\|\bar U_{t_k}\|
\le {\bar M_*}.
\]
Arguing as in the proof of \Cref{thm:stability1},
\(
\|\bar D_0\|_2\le\frac{1}{\bar m_*}\|\bar Z_0\|_2.
\)
Combining this with \eqref{eq:numVUdecay}, we obtain
\[
\|\bar Z_{t_k}\|^2_2
\le \frac{\bar M_*^2}{\bar m^2_*}e^{-\frac{b\,}{1+b\Delta t}k\Delta t} \|\bar{Z}_0\|_2^2.
\]
Thus, the desired estimate follows with \(C:={{\bar M_*^2}}/{\bar m^2_*}.\)
\end{proof}

\begin{remark}
    The above theorem shows that for any step size $\Delta t$, the Euler scheme \eqref{eq:freeCIRfem} exponentially stable with decay rate $\lambda_2 = \frac{b}{1+b\Delta t}$, which converges to the continuous rate $b$ in \Cref{thm:stability1} as the step size tends to zero.
\end{remark}
\section{Experimental results}
\label{sec:experiments}
In this section, we present several numerical experiments to illustrate the convergence, long-time stability, and positivity preservation of the proposed drift-implicit square-root Euler scheme. Following  \cite{numfsde1,numfsde2}, we approximate the corresponding free processes using self-adjoint random matrices of sufficiently large finite dimension $N$. We emphasize, however, that our theoretical error analysis concerns only time discretization and does not quantify the finite-N random matrix approximation error.

Throughout the numerical experiments, we approximate \(\varphi\) by the normalized matrix trace \(\operatorname{tr}_N(\cdot):=\frac{1}{N}\operatorname{tr}(\cdot),\) the free Brownian increments by
\[
\Delta W_k^{(N)}
=
\sqrt{\frac{\Delta t}{2N}}
\left(G_k^{(N)}+\bigl(G_k^{(N)}\bigr)^*\right),
\quad k=0,\ldots,L-1,
\]
where $G_k^{(N)}$ are independent $N\times N$ random matrices with independent standard Gaussian entries. Unless otherwise stated, the parameters in \eqref{eq:freeCIRmain} are \(a=2,\ b=4,\ \sigma=1,\ X_0^{(N)}=V_0^{(N)}=I_N.\)

By \Cref{prop:wellposed}, our scheme admits the closed-form update
\[
V_{k+1}^{(N)}=\frac{B_k+\sqrt{B_k^2+\xi\Delta t\left(2+b\Delta t\right)I_N}}{2+b\Delta t},\quad B_k=V_k^{(N)}+\frac{\sigma}{2}\Delta W_{k}^{(N)},
\]
and the corresponding approximation for \eqref{eq:freeCIRmain} is \(X_k^{(N)}=\big(V_k^{(N)}\big)^2.\)

\begin{figure}[!ht]
    \centering
    \begin{subfigure}[b]{0.35\textwidth}
        \includegraphics[width=\textwidth]{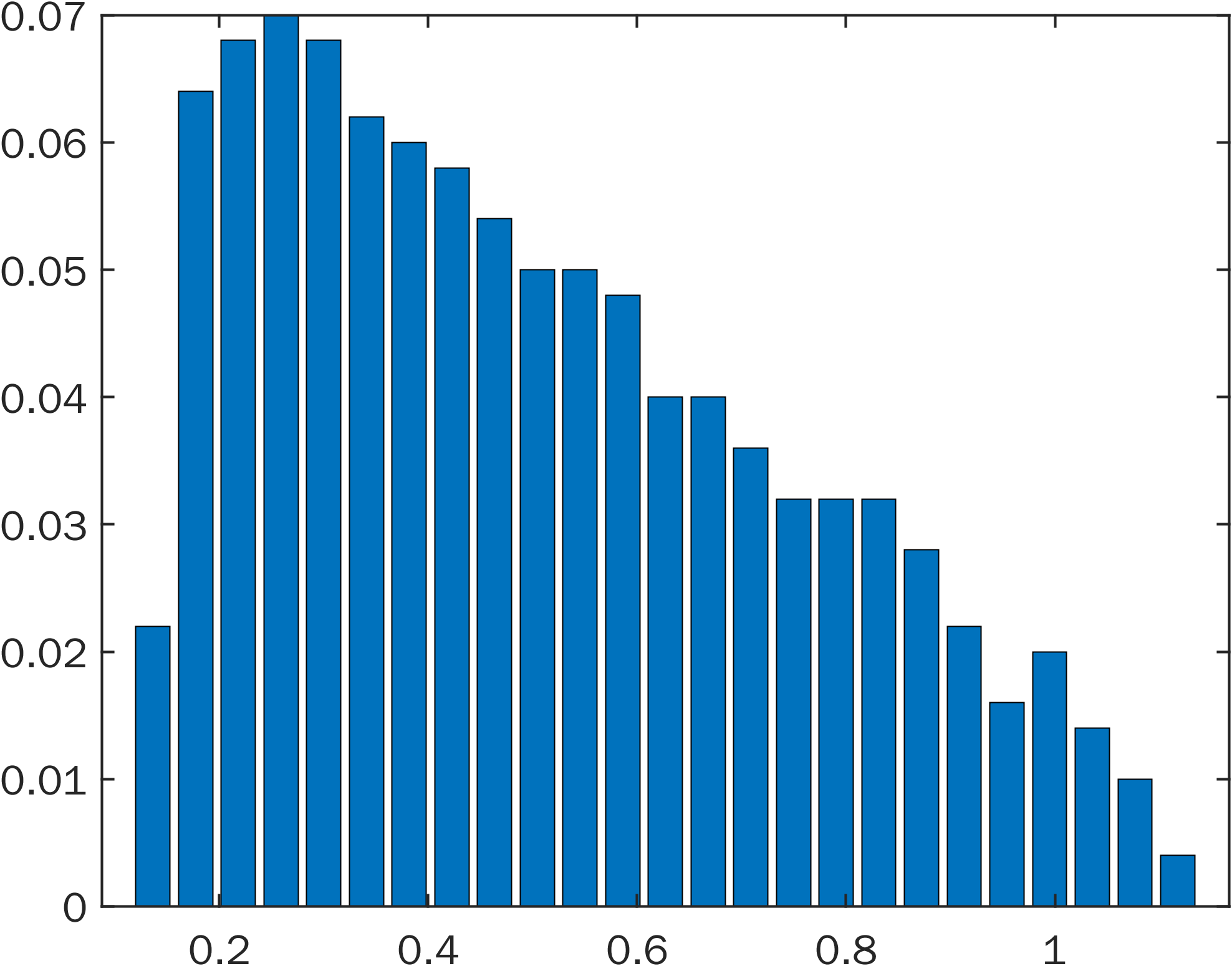}
        \caption{L=300.}
        \label{fig:sub1}
    \end{subfigure}
    \begin{subfigure}[b]{0.35\textwidth}
        \includegraphics[width=\textwidth]{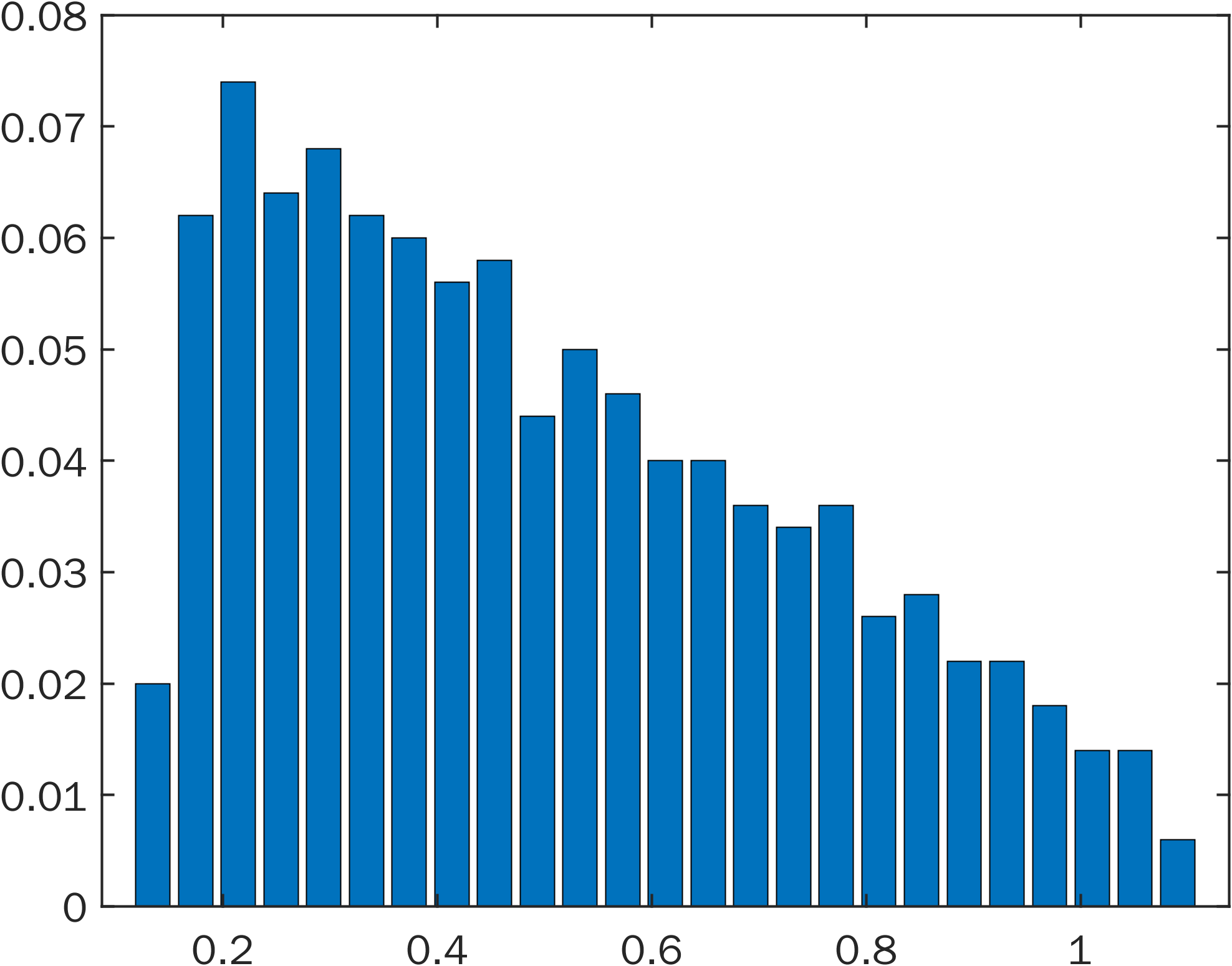}
        \caption{L=4000.}
        \label{fig:sub2}
    \end{subfigure}
    \caption{Empirical spectral distributions of $X_L^{(N)}$ with different time steps.}
    \label{fig:total}
\end{figure}

 It is shown in \Cref{fig:total} that the empirical spectral distributions of \(X_L^{(N)}\) with \(L=300\), \(L=4000\), \(N=500\) and \(T=1\). Applying $\varphi$ to \eqref{eq:freeCIRmain} gives
\[
\varphi(X_T)=\frac12+\frac12e^{-4}\approx0.50915782.
\]
\Cref{tab:first_moment} shows the expectation of the spectrum distribution of the numerical solution tends to the expectation of the exact solution as $N$ increases.

\begin{table}[!ht]
\centering
\renewcommand{\arraystretch}{1}
\setlength{\tabcolsep}{18pt}
\caption{Monte Carlo mean of $\operatorname{tr}_N(X^{(N)}_L)$ and its absolute error with respect to $\varphi(X_T)$ for different matrix dimensions $N$.}
\label{tab:first_moment}
    \begin{tabular}{c|c|c}\hline
        \rule[-1.2ex]{0pt}{4ex}$N$ 
        & $\frac{1}{M}\sum_{m=1}^M\operatorname{tr}_N(X^{(N,m)}_L)$
        & Absolute error \\  \hline
        50  & 0.51021257 & $1.0548\times10^{-3}$ \\
        200 & 0.50933525 & $1.7743\times10^{-4}$ \\  
        500 & 0.50922616 & $6.8339\times10^{-5}$ \\  \hline
    \end{tabular}
\end{table}

\begin{figure}[htbp]
    \centering
    \begin{subfigure}[b]{0.45\textwidth}
        \includegraphics[width=\textwidth]{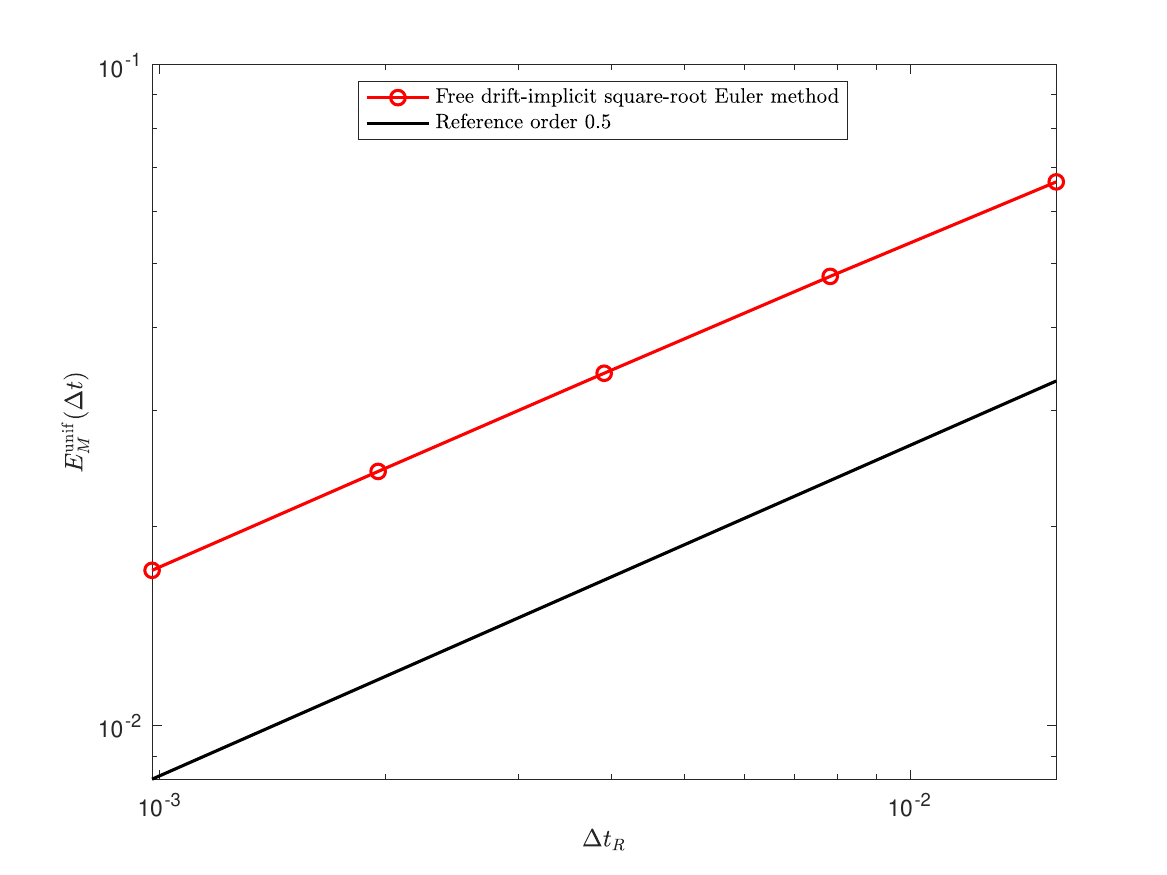}
    \caption{Uniform-in-time  error}
    \end{subfigure}
    \begin{subfigure}[b]{0.45\textwidth}
        \includegraphics[width=\textwidth]{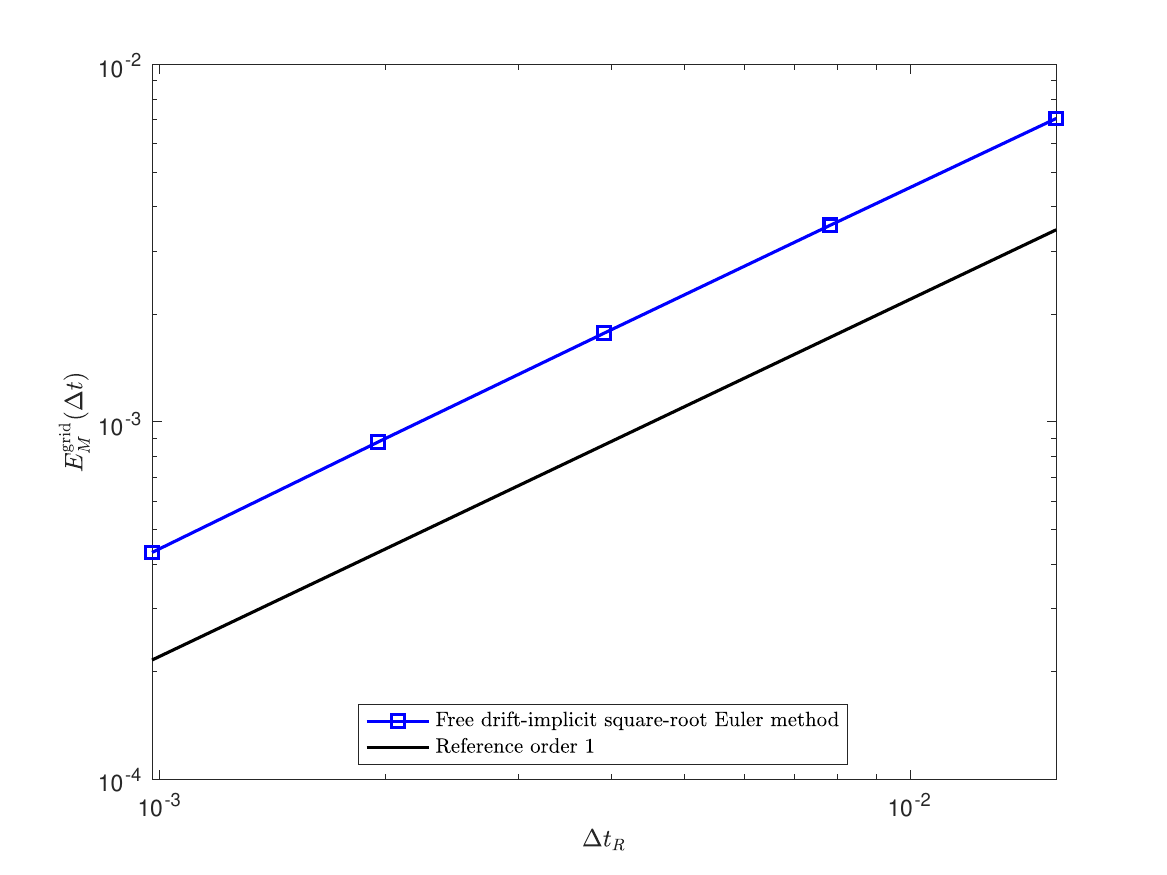}
    \caption{Grid-points error}
    \end{subfigure}
    \caption{Strong convergence of the free drift-implicit square-root Euler scheme.}
    \label{fig:strong_convergence}
\end{figure}

For strong convergence, we use a reference solution with \(\Delta t_{\mathrm{ref}}=2^{-14}\) and coarse step sizes
\(
\Delta t\in\{2^{-10},2^{-9},2^{-8},2^{-7},2^{-6}\}.
\)
The coarse and reference solutions are driven by the same matrix Brownian path, with \(\Delta W_k^{(N,\Delta t)}=\sum_{j=kR}^{(k+1)R-1}\Delta W_j^{(N,\Delta t_{\mathrm{ref}})}\) and \(R=\Delta t/\Delta t_{\mathrm{ref}}.\)
We take $N=100$ and use $M=100$ independent realizations.  For each coarse step size $\Delta t$, the strong error is defined by
\[
E_M(\Delta t;\mathcal T)
:=
\max_{t\in\mathcal T}
\left(
\frac{1}{M}\sum_{m=1}^M
\bigl\|
\bar X_{\Delta t}^{(m)}(t)
-
X_{\mathrm{ref}}^{(m)}(t)
\bigr\|_{2,N}^2
\right)^{\frac12},
\]
where \(\|\cdot \|^2_{2,N}:=\operatorname{tr}_N(|\cdot|^2).\)
We take $\mathcal T=\{t_k\}_{k=1}^{L}$ for the grid-point error and $\mathcal T=\{t_k+\Delta t/2\}_{k=0}^{L-1}$ for the midpoint approximation of the
uniform-in-time error.
It is shown in \Cref{fig:strong_convergence} that the uniform-in-time error and the grid-point error decrease at slopes approaching  \(1/2\) and \(1,\) respectively, consistent with  \Cref{thm:convergence}.

\begin{figure}[htbp]
    \centering
        \begin{subfigure}[b]{0.45\textwidth}
        \includegraphics[width=\textwidth]{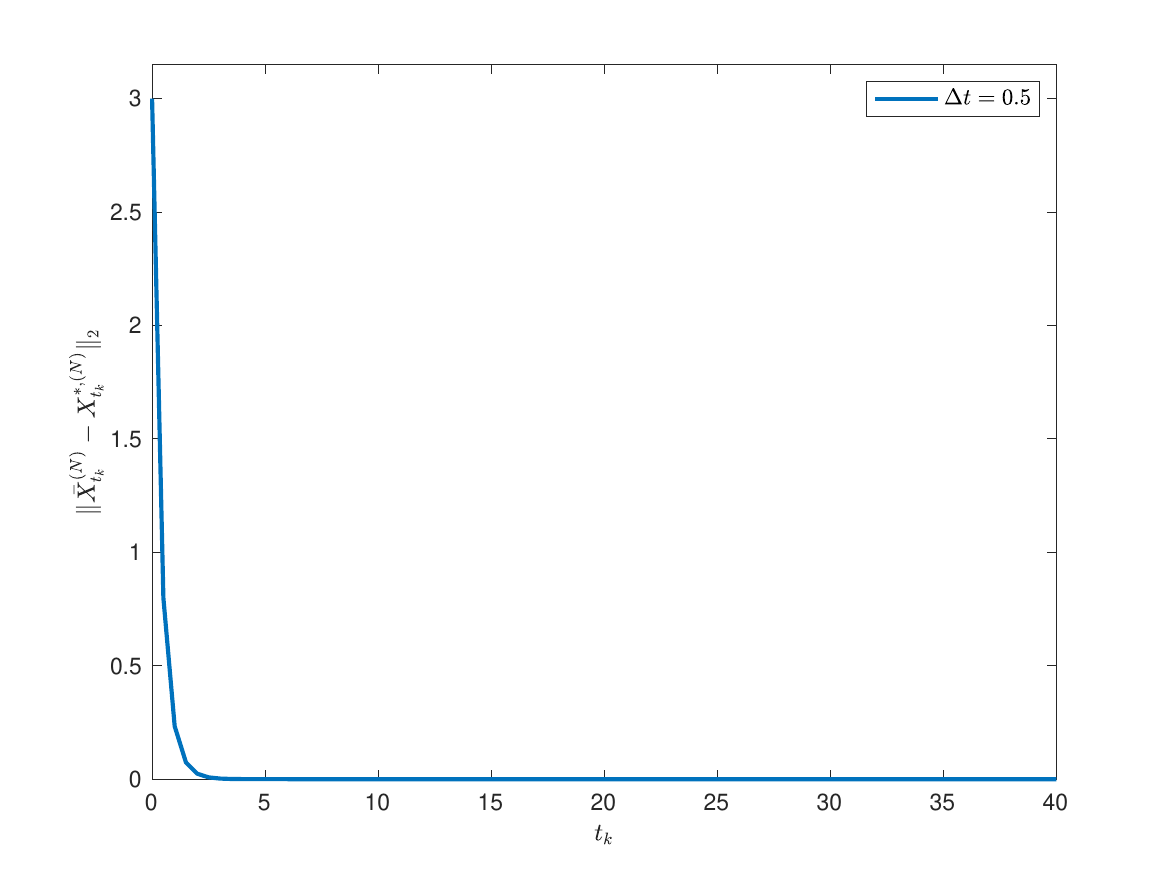}
        \caption{$\Delta t=\frac{1}{2}$}
        \label{fig:sub3}
    \end{subfigure}  
    \begin{subfigure}[b]{0.45\textwidth}
        \includegraphics[width=\textwidth]{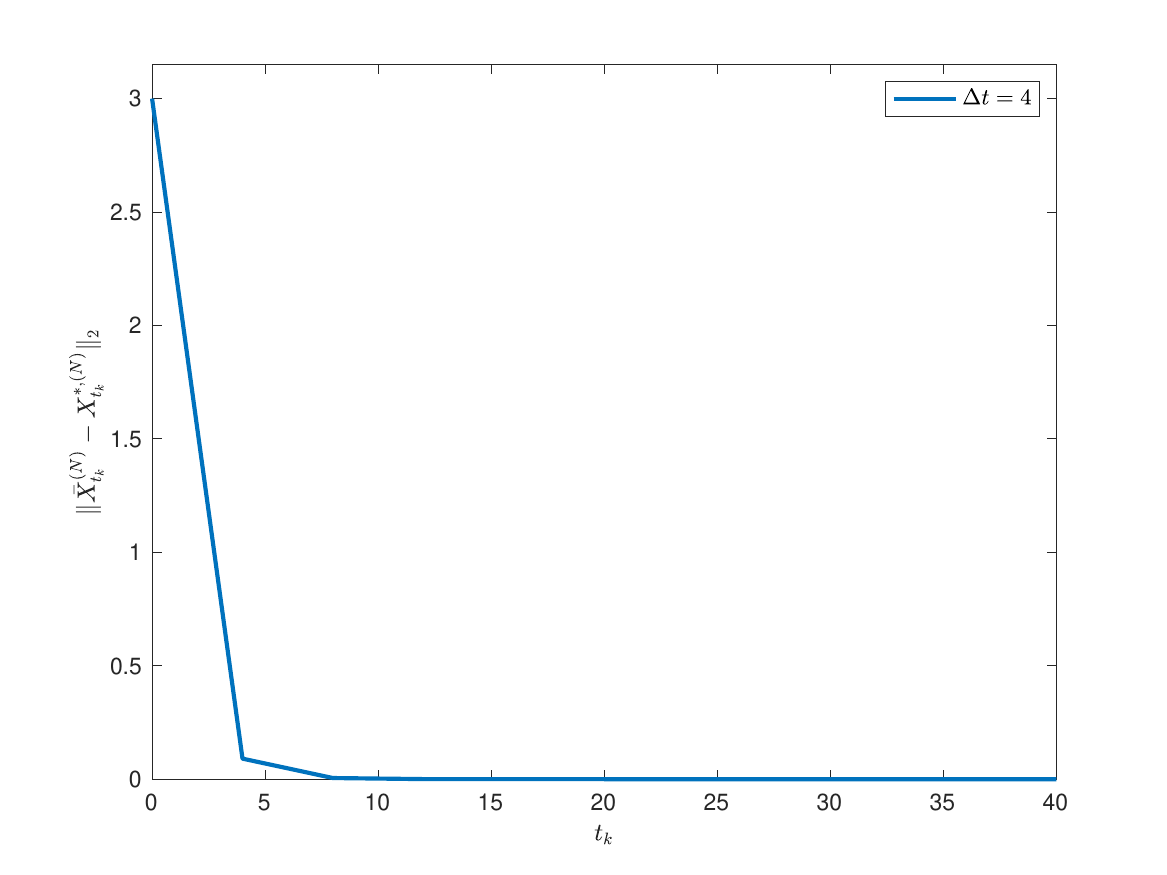}
        \caption{$\Delta t=4$}
        \label{fig:sub4}
    \end{subfigure}
    \caption{Exponential stability of the proposed scheme for different step sizes.}
    \label{fig:stability}
\end{figure}

For stability, we compare solutions starting from $X_0^{(N)}=I_N$ and $\widetilde X_0^{(N)}=4I_N$, driven by the same matrix Brownian increments, with \(T=40\) and \(\Delta t=1/2,\ 4\). \Cref{fig:stability} shows decay of $\bigl\|X_k^{(N)}-\widetilde X_k^{(N)}\bigr\|_{2,N}$
for both step sizes, consistent with the unconditional exponential stability proved in \Cref{thm:numstability}.

\begin{figure}[htbp]
    \centering
     \includegraphics[width=0.5\textwidth]{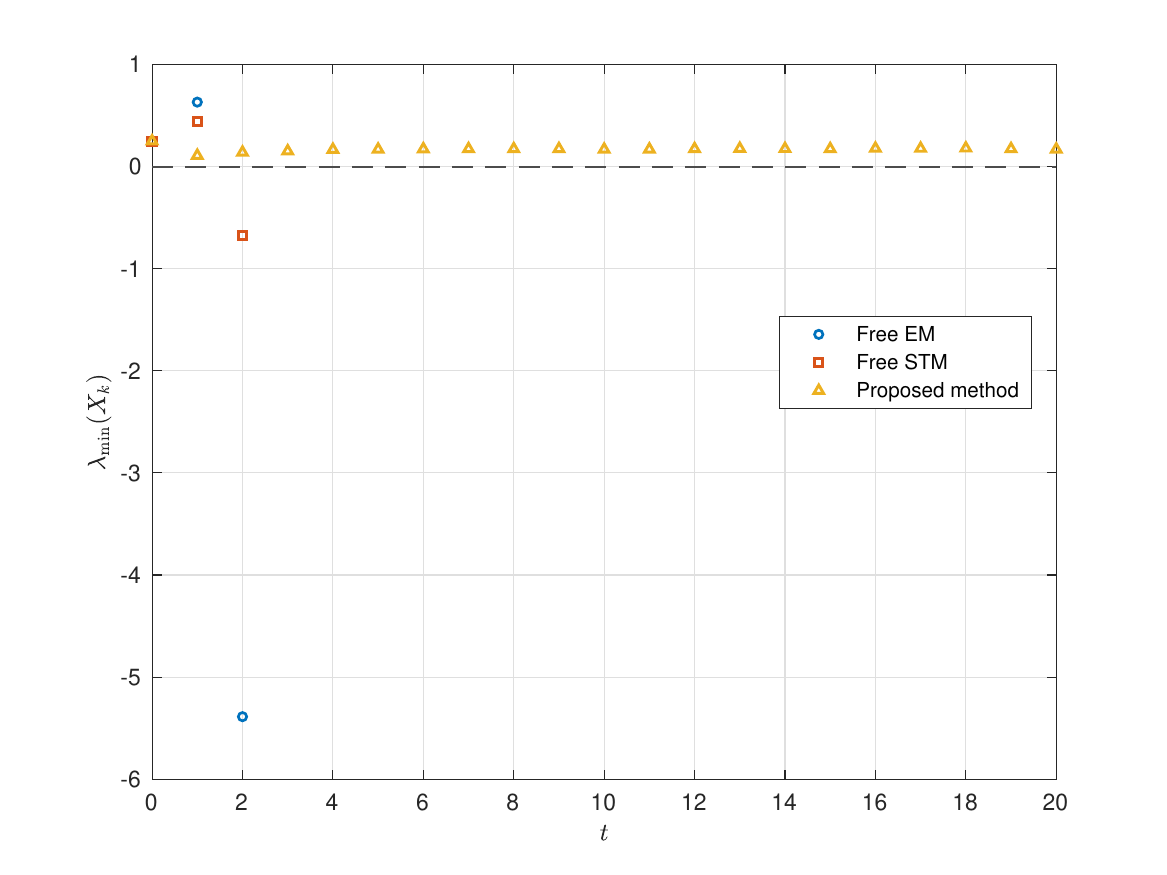}
    \caption{Comparison of the positivity-preserving behavior of the free EM, free STM, and the proposed free drift-implicit square-root Euler method.}
    \label{fig:method_comparison}
\end{figure}

Finally, we compare the positivity preservation of the proposed scheme with that of the free Euler-Maruyama method \cite{numfsde1} and the free stochastic theta method \cite{numfsde2} applied to \eqref{eq:freeCIRmain}. We take
\[
a=3,\quad b=1,\quad \sigma=2.4,\quad
N=500,\quad T=20,\quad \Delta t=1,
\]
and \(\theta=1\) for the theta method. \Cref{fig:method_comparison} shows that the proposed scheme maintains strictly positive eigenvalues, illustrating the positivity preservation established in \Cref{prop:wellposed}, whereas the other two methods failed.

\end{document}